\documentclass[11pt,leqno]{article}

\usepackage{amsmath,amsthm,amssymb}
\usepackage{authblk,subfigure,booktabs}
\usepackage{enumerate}
\usepackage{graphicx}
\usepackage{color}

\usepackage{mathabx,bm}
\usepackage{algorithm}
\usepackage[noend]{algpseudocode}
\usepackage{url}
\usepackage{enumitem}

\usepackage{geometry}
\usepackage{listings}

\usepackage{apptools}

\theoremstyle{plain}
\newtheorem{thm}{Theorem}[section]

\newtheorem{lem}[thm]{Lemma}

\theoremstyle{remark}
\newtheorem{rem}[thm]{Remark}

\theoremstyle{definition}

\makeatletter

\numberwithin{equation}{section}
\numberwithin{figure}{section}
\numberwithin{table}{section}

\title{\Large\bf Convergence of differentiable non-monotone schemes \\ for fully nonlinear parabolic equations}

\author{Yumiharu Nakano\thanks{E-mail: nakano@comp.isct.ac.jp}}
\affil{Department of Mathematical and Computing Science \protect \\ Institute of Science Tokyo}

\date{\today}

\begin{document}

\maketitle

\begin{abstract}
We develop a convergence theory for non-monotone approximation schemes for
fully nonlinear parabolic partial differential equations.
Modern computational methods such as kernel-based collocation, spectral methods,
physics-informed neural networks, and deep Galerkin methods are typically non-monotone,
since they produce smooth approximate solutions and compute spatial derivatives
directly from gradients of the chosen ansatz.
Such schemes therefore lie outside the scope of the classical Barles and Souganidis
convergence theory.
We introduce an abstract framework that replaces strict monotonicity by two pointwise
consistency conditions, on the PDE residual and on the terminal mismatch, both directly
verifiable for a smooth approximating sequence.
The technical key is a max-min representation of the nonlinearity that converts a
vanishing classical residual into the viscosity subsolution and supersolution
inequalities, and so dispenses with monotonicity in the abstract argument.
The framework yields qualitative convergence under standard hypotheses,
together with a quantitative error bound for Hamilton-Jacobi-Bellman equations on
an unbounded spatial domain, in which the residual on an expanding truncation cylinder
is balanced against an exponentially decaying tail term coming from the controlled
stochastic differential equation underlying the value function.
As a concrete realization, we analyze kernel-based collocation with Wendland radial basis
functions, and present numerical experiments on a benchmark Hamilton-Jacobi-Bellman
problem in one and two spatial dimensions that confirm the predicted convergence behaviour.

\begin{flushleft}
\end{flushleft}
\begin{flushleft}
\end{flushleft}
\end{abstract}





\section{Introduction}\label{sec:1}

This paper develops a convergence theory for non-monotone approximation schemes
applied to terminal value problems of fully nonlinear parabolic partial differential equations:
\begin{equation}
\label{eq:1.1}
\left\{
\begin{split}
 &-\partial_t v + F(t,x,v(t,x),D_xv(t,x),D^2_{xx}v(t,x))=0, \quad (t,x)\in
 [0,T)\times \mathbb{R}^d, \\
 &v(T,x)=g(x), \quad x\in\mathbb{R}^d,
\end{split}
\right.
\end{equation}
where $F:[0,T]\times\mathbb{R}^d\times\mathbb{R}\times\mathbb{R}^d\times\mathbb{S}^d\to\mathbb{R}$
is degenerate elliptic, and $\mathbb{S}^d$ denotes the space of symmetric $d\times d$ real matrices.
Throughout, $\partial_t$ denotes the partial derivative in $t$, and $D_x$, $D^2_{xx}$ denote
the gradient and Hessian in $x$, respectively.

Problem \eqref{eq:1.1} arises naturally in stochastic optimal control: when $F$ is of
Hamilton--Jacobi--Bellman (HJB) type, the solution $v$ is the value function of a stochastic
control problem and an optimal control policy can be recovered from $v$.
Owing to the possible degeneracy and nonlinearity of $F$, classical smooth solutions cannot
be expected in general, and \eqref{eq:1.1} is interpreted in the viscosity sense.
Under standard structural hypotheses including degenerate ellipticity, \eqref{eq:1.1}
admits a unique continuous viscosity solution; see Fleming and Soner~\cite{fle-son:2006},
Pham~\cite{pha:2009}, Kohn and Serfaty~\cite{koh-ser:2010}.

\medskip
\noindent\textbf{Smooth, non-monotone schemes.}
For a long time, the standard tool for proving convergence of numerical methods for
\eqref{eq:1.1} has been the abstract framework of Barles and Souganidis~\cite{bar-sou:1991},
which requires the scheme to be monotone, stable, and consistent.
Classical schemes designed to satisfy this triple include the finite difference methods of
Kushner and Dupuis~\cite{kus-dup:2001} and Bonnans and Zidani~\cite{bon-zid:2003},
the semi-Lagrangian methods of Camilli and Falcone~\cite{cam-fal:1995} and
Debrabant and Jakobsen~\cite{deb-jak:2012}, and the probabilistic methods of
Pag{\`e}s et al.~\cite{pag-pha-pri:2004}, Fahim et al.~\cite{fah-tou-war:2011},
Guo et al.~\cite{guo-etal:2015}, and Nakano~\cite{nak:2014b}.
A more recent generation of computational methods, however, parameterizes the approximate
solution by a smooth ansatz and computes spatial derivatives directly from gradients of that
ansatz: kernel-based collocation methods of Kansa~\cite{kan:1990b} and Nakano~\cite{nak:2017},
spectral and Galerkin methods, deep Galerkin methods of Sirignano and Spiliopoulos~\cite{sir-spi:2018},
and physics-informed neural networks of Raissi, Perdikaris and Karniadakis~\cite{rai-per-kar:2019}.
These smooth-ansatz schemes are attractive because they offer mesh-free, dimension-flexible,
and easily implementable approximations; their derivatives are produced by automatic
differentiation or by explicit gradient formulas of the basis. They are, however, typically
\emph{non-monotone}: smooth bases inevitably mix positive and negative weights when
representing the second-derivative operator, so the strict monotonicity condition of
Barles--Souganidis is violated.

For quasilinear parabolic equations, Sirignano and Spiliopoulos~\cite{sir-spi:2018} obtained
convergence of a deep Galerkin scheme to the viscosity solution under the assumption that the
$L^2$-loss can be made arbitrarily small along a sequence of approximators with uniformly
bounded $C^{1,2}$-norm; the existence of such a sequence and the actual behaviour of the
optimization remain implicit, and the analysis does not extend to fully nonlinear $F$.
For Kolmogorov-type linear or semilinear equations, error analyses of physics-informed neural
networks have been developed by Mishra and Molinaro~\cite{mis-mol:2022} and
De Ryck and Mishra~\cite{drm:2024}, but again do not address the fully nonlinear viscosity
solution setting. To the author's knowledge, no general convergence theory in the
Barles--Souganidis spirit is available for non-monotone smooth schemes applied to fully
nonlinear parabolic equations on the unbounded spatial domain~$\mathbb{R}^d$.

\medskip
\noindent\textbf{Contributions.}
The aim of this paper is to fill that gap, with the following four contributions.

\smallskip
\noindent\emph{(i) An abstract convergence framework based on pointwise consistency.}
We introduce an abstract scheme determined by a sequence $\{v_n\}$ of smooth candidate solutions
on $[0,T]\times\mathbb{R}^d$ and impose two consistency conditions, \rm{(B1)} on the PDE residual
and \rm{(B2)} on the terminal mismatch, both pointwise (Section~\ref{sec:2.1}).
Theorem~\ref{thm:2.1} establishes that under these conditions and the standard hypotheses
on $F$, the sequence $\{v_n\}$ converges, uniformly on every compact subset of
$[0,T]\times\mathbb{R}^d$, to the unique continuous viscosity solution of \eqref{eq:1.1}.
The framework dispenses with monotonicity entirely.

\smallskip
\noindent\emph{(ii) A max-min representation of $F$ as the technical key.}
The proof of Theorem~\ref{thm:2.1} rests on a max-min representation of the nonlinearity
$F$ established in our earlier work~\cite{nak:2017}: at any $C^3$ test function $\varphi$,
the residual $\varphi - h F(\cdot;\varphi)$ is approximated, up to $O(h^{1+\beta})$, by a
sup--inf over a parametric drift--diffusion functional. This representation converts a small
classical residual into the viscosity sub/supersolution inequality at the contact point with
the test function, even when the scheme is non-monotone.
Lemma~\ref{lem:2.2} re-states this tool in the form needed here.

\smallskip
\noindent\emph{(iii) A quantitative error bound for Hamilton--Jacobi--Bellman equations.}
For HJB equations on $[0,T]\times\mathbb{R}^d$, where $F$ is the supremum over a control set
of a controlled drift--diffusion generator, Theorem~\ref{thm:2.2} provides a quantitative
estimate
\[
 \sup_Q|v - v_n| \le (1+T)\eta_n + C\lambda_n e^{-R_n^2/(8 C_1^2 T)},
\]
where $\eta_n$ is a sup-norm bound on the residual and terminal mismatch on a truncation
cylinder of radius $R_n$, the constant $\lambda_n$ controls the $C^{1,2}$-norm of $v_n$,
and the exponential tail comes from a Gaussian concentration estimate for the controlled SDE.
The bound is uniform on any compact $Q$ and accommodates both fixed and growing
$\lambda_n,R_n$ chosen by the user.

\smallskip
\noindent\emph{(iv) Realization with kernel-based collocation.}
To illustrate the abstract framework, Section~\ref{sec:3} treats kernel-based collocation
with a positive-definite radial basis function $\Phi$ whose Fourier transform decays like
$(1+|\xi|^2)^{-\tau}$. Approximate solutions take the form $v_n(y)=\sum_j\theta_j\Phi(y-y^j)$
on a fill-distance grid, and we obtain qualitative convergence for general $F$
(Theorem~\ref{thm:3.1}) and a fully explicit error rate for HJB equations
(Theorem~\ref{thm:3.3}). The required $C^{1,2}$ control of $v_n$ is supplied by the native space
norm constraint $\theta^{\mathsf T}K\theta \le \lambda$ via Sobolev embedding.
We close with numerical experiments using Wendland radial basis functions on a benchmark
HJB equation that confirm the predicted convergence behaviour and assess the computational cost.

\medskip
\noindent\textbf{Organization.}
Section~\ref{sec:2} develops the abstract framework: Section~\ref{sec:2.1} treats general $F$
(qualitative convergence, Theorem~\ref{thm:2.1}) and Section~\ref{sec:2.2} treats the HJB
case (quantitative bound, Theorem~\ref{thm:2.2}).
Section~\ref{sec:3} applies the framework to kernel-based collocation; Section~\ref{sec:4}
contains the numerical experiments.

{\it Notation.}
Throughout this paper,
$a^{\mathsf{T}}$ denotes the transpose of a vector or matrix $a$.
For $a=(a_i)\in\mathbb{R}^{\ell}$ we write
$|a|=(\sum_{i=1}^{\ell}a_{i}^2)^{1/2}$,
$|a|_{\infty}=\max_{i=1,\ldots,\ell}|a_i|$, and $|a|_1=\sum_{i=1}^{\ell}|a_i|$.
For $a=(a_{ij})\in\mathbb{R}^{m\times k}$ we write
$|a| = (\sum_{i=1}^m\sum_{j=1}^k a_{ij}^2)^{1/2}$.

For a multi-index $\alpha=(\alpha_1,\ldots,\alpha_d)$ of nonnegative integers,
the differential operator $D^{\alpha}$ is defined by
\begin{equation*}
 D^{\alpha}f(x)=\frac{\partial^{|\alpha|_1}}
  {\partial x_1^{\alpha_1}\cdots \partial x_d^{\alpha_d}}f(x), \quad x\in\mathbb{R}^d.
\end{equation*}
We write $D^{\alpha}_xf$ for $D^{\alpha}f$ when we wish to emphasize
differentiation in $x$.
We also write $\partial_t^mz(t,x):= (\partial^m z/\partial t^m)(t,x)$.
For the time-space variable, we use the extended multi-index
$\alpha=(\alpha_0,\alpha_1,\ldots,\alpha_d)$ and set
\[
 D^{\alpha}z(t,x)=\frac{\partial^{|\alpha|_1}}
  {\partial_t^{\alpha_0}\partial x_1^{\alpha_1}\cdots \partial x_d^{\alpha_d}}z(t,x).
\]
For an open set $\mathcal{O}\subset\mathbb{R}^d$ or $\mathbb{R}^{d+1}$,
$C^{\kappa}(\mathcal{O})$ denotes the space of real-valued functions
with continuous derivatives up to order $\kappa\in\mathbb{N}\cup\{0\}$,
with norm
\[
 \|f\|_{C^{\kappa}(\mathcal{O})}=\sum_{|\beta|_1\le\kappa}\sup_{x\in\mathcal{O}}|D^{\beta}f(x)|.
\]
We write $C(\mathcal{O})=C^0(\mathcal{O})$.
For $\nu,\kappa\in\mathbb{N}\cup\{0\}$ and $t_0<t_1$,
$C^{\nu,\kappa}([t_0,t_1]\times O)$ denotes the space of functions that are
$\nu$-times continuously differentiable in $t$ and
$\kappa$-times continuously differentiable in $x$, with norm
\[
 \|z\|_{C^{\nu,\kappa}([t_0,t_1]\times O)}
 = \|z\|_{C^0([t_0,t_1]\times O)}
   + \sum_{m=1}^{\nu}\|\partial_t^m z\|_{C^0([t_0,t_1]\times O)}
   + \sum_{1\le |\alpha|_1\le\kappa}\|D^{\alpha}_xz\|_{C^0([t_0,t_1]\times O)}.
\]
$C^{\infty}_0(\mathbb{R}^m)$ denotes the space of smooth functions
with compact support in $\mathbb{R}^m$.
For $(t,x)\in\mathbb{R}\times\mathbb{R}^d$, we sometimes write $y=(t,x)\in\mathbb{R}^{d+1}$
and $f(t,x)=f(y)$ when convenient.

\section{Abstract convergence framework}\label{sec:2}

\subsection{General equations}\label{sec:2.1}

 We study the convergence of the approximation method described in Section \ref{sec:3}
under the conditions where \eqref{eq:1.1} admits a unique viscosity solution.  
To this end, first we recall the notion of the viscosity solution and 
describe our standing assumptions for \eqref{eq:1.1}. 

An $\mathbb{R}$-valued, upper-semicontinuous function 
$u$ on $[0,T]\times\mathbb{R}^d$ is said to be a viscosity subsolution of 
(\ref{eq:1.1}) if the following two conditions hold: 
\begin{itemize}
 \item[(i)] for every $(t,x)\in [0,T)\times\mathbb{R}^d$ and 
every $\varphi\in C^{1,2}([0,T]\times\mathbb{R}^d)$ such that 
$0=(u-\varphi)(t,x)=\max_{(s,y)\in [0,T)\times\mathbb{R}^d}(u-\varphi)(s,y)$ we have 
\begin{equation*}
 -\partial_t\varphi(t,x)+F(t,x,u(t,x),D_x\varphi(t,x),
  D_{xx}^2\varphi(t,x))\le 0; 
\end{equation*}
 \item[(ii)] $u(T,x)\le g(x)$, $x\in\mathbb{R}^d$. 
\end{itemize}
Similarly, an $\mathbb{R}$-valued, lower-semicontinuous function 
$u$ on $[0,T]\times\mathbb{R}^d$ is said to be a viscosity supersolution of 
(\ref{eq:1.1}) if the following two conditions hold: 
\begin{itemize}
 \item[(i)] for every $(t,x)\in [0,T)\times\mathbb{R}^d$ and 
every $\varphi\in C^{1,2}([0,T]\times\mathbb{R}^d)$ such that 
$0=(u-\varphi)(t,x)=\min_{(s,y)\in [0,T]\times\mathbb{R}^d}(u-\varphi)(s,y)$ we have 
\begin{equation*}
 -\partial_t\varphi(t,x)+F(t,x,u(t,x),D_x\varphi(t,x),
   D_{xx}^2\varphi(t,x))\ge 0;  
\end{equation*}
 \item[(ii)] $u(T,x)\ge g(x)$, $x\in\mathbb{R}^d$. 
\end{itemize}
We say that $u$ is a viscosity solution of (\ref{eq:1.1}) if 
it is both a viscosity subsolution and a viscosity supersolution 
of (\ref{eq:1.1}). 

We consider the terminal value problem (\ref{eq:1.1}) 
under the following assumptions: 
there exists a positive constant $C_0$ such that the following are satisfied: 
\begin{enumerate}
\item[(A1)] For $t\in [0,T]$, $x\in\mathbb{R}^d$, $z\in\mathbb{R}$, 
$p,p^{\prime}\in\mathbb{R}^d$, and $X,X^{\prime}\in\mathbb{S}^d$, 
  \[
   |F(t,x,z,p,X)- F(t,x,z,p^{\prime},X^{\prime})| \le C_0(|p-p^{\prime}|+|X-X^{\prime}|).    
  \]
\item[(A2)] For $t\in [0,T]$, $x\in\mathbb{R}^d$, $z\in\mathbb{R}$, 
$p\in\mathbb{R}^d$, and $X\in\mathbb{S}^d$, 
\begin{equation*}
|F(t,x,z,p,X)|\le C_0(1+|z|+|p|+|X|). 
\end{equation*}
 \item[(A3)] For $t\in [0,T]$, $x\in\mathbb{R}^d$, $z\in\mathbb{R}$, 
$p\in\mathbb{R}^d$, and $X,X^{\prime}\in\mathbb{S}^d$ with 
$X\ge X^{\prime}$, 
  \begin{equation*}
   F(t,x,z,p,X)\le F(t,x,z,p,X^{\prime}).  
  \end{equation*}
\item[(A4)] The function $g$ is Lipschitz continuous and bounded on $\mathbb{R}^d$. 
\end{enumerate}

We assume that the following comparison principle holds: 
\begin{enumerate}
\item[(A5)] For every bounded, upper-semicontinuous viscosity subsolution $u$ 
of (\ref{eq:1.1}) 
and bounded lower-semicontinuous viscosity supersolution $w$ of 
(\ref{eq:1.1}), we have 
\begin{equation*}
 u(t,x)\le w(t,x), 
 \quad (t,x)\in [0,T]\times\mathbb{R}^d. 
\end{equation*}
\end{enumerate}

The conditions (A1)--(A5) ensure that
there exists a unique continuous viscosity solution $v$ of (\ref{eq:1.1}); see \cite{koh-ser:2010}.

Let $\{v_n\}_{n=1}^{\infty}$ be a given sequence of functions such that $\|v_n\|_{C^{2,3}([0,T)\times O)}<\infty$ for any 
bounded open $O\subset\mathbb{R}^d$ and for any $n\in\mathbb{N}$. 
Then we make the following conditions on our scheme: 
\begin{enumerate}
\item[(B1)] {\it Consistency I}: for any $(t,x)\in [0,T)\times\mathbb{R}^d$, 
\[
 \lim_{{(s,y)\to (t,x)}\atop{n\to\infty}}\left|-\partial_sv_n(s,y) + F(s,y; v_n(s,\cdot))\right|=0, 
\]
where for any $C^2$-function $\varphi$ on $\mathbb{R}^d$,  
\[
 F(t,x;\varphi) = F(t,x,\varphi(x), D_x\varphi(x),D^2_{xx}\varphi(x)), 
  \quad x\in\mathbb{R}^d.
\]
\item[(B2)] {\it Consistency II}: for $x\in\mathbb{R}^d$, 
\[
 \lim_{{(s,y)\to (T,x)}\atop{n\to\infty}}v_n(s,y) = g(x). 
\]
\end{enumerate}

Now we are ready to state one of our main results, which establishes the convergence of
the smooth abstract methods.
\begin{thm}
\label{thm:2.1}
Suppose that \rm{(A1)}--\rm{(A5)}, \rm{(B1)}, and \rm{(B2)} hold. 
Then  
\begin{equation*}
 v_n(s,x)\to v(t,x),   
\end{equation*}
as $s\to t$ and $n\to\infty$ uniformly on any compact subset of $\mathbb{R}^d$. 
\end{thm}

The rest of this section is devoted to the proof of Theorem \ref{thm:2.1}. 
In what follows, by $C$ we denote positive constants that may vary from line to line 
and that are independent of $n$ and $(t,x)\in [0,T]\times\mathbb{R}^d$.

For $h>0$, $c>0$ and $\kappa>0$ define 
\begin{equation*}
 \mathcal{D}_{h, c}=\left\{(p,\Gamma)\in\mathbb{R}^d\times 
  \mathbb{S}^d: |p|, |\Gamma|\le c h^{-q}\right\}, \quad 
 \mathcal{X}_{h,\kappa}=\left\{w\in\mathbb{R}^d: |w|\le h^{-\kappa}\right\}. 
\end{equation*}
Our key tool for proving Theorem \ref{thm:2.1} is the following result, which is proved in \cite{nak:2017}.  
\begin{lem}[Lemma 3.12 in \cite{nak:2017}]
\label{lem:2.2}
Suppose that \rm{(A1)}--\rm{(A4)} hold. Let $O$ be a bounded open subset of $\mathbb{R}^d$. 
Then for any $q\in (0,1/5)$ and open ball $U$ compactly included in $O$ 
there exist $h_0\in (0,1]$, $\beta\in (0,\infty)$, 
$\kappa\in (q/2,\infty)$ such that 
for $(t,x,z)\in [0,T]\times U\times \mathbb{R}$,  
$\varphi\in C^3(O)$ with 
$\|\varphi\|_{C^3(O)}
\le ch^{-q}$, and $h\le h_0$, 
\begin{align*}
 &\bigg|\varphi(x)-h F(t,x,z,D\varphi(x),D^2\varphi(x)) \\
 & \quad -\sup_{(p,\Gamma)\in\mathcal{D}_{h,c}}\inf_{w\in\mathcal{X}_{h,\kappa}}
  \left[\varphi(x+\sqrt{h}w) 
 - \sqrt{h}w^{\mathsf{T}}p-\frac{h}{2}w^{\mathsf{T}}\Gamma w  
 - h F(t,x,z,p,\Gamma)\right]\bigg| \le C h^{1+\beta}. 
\end{align*}
\end{lem}

Under our consistency conditions, the function $v_n$ is actually locally bounded uniformly in $n$. 
\begin{lem}
\label{lem:2.3}
Under the assumptions imposed in Theorem $\ref{thm:2.1}$, 
for any open ball $U\subset\mathbb{R}^d$ 
there exist $n_0\in\mathbb{N}$ such that 
\[
 \sup_{n\ge 1}\|v_n\|_{C([0,T]\times U)}<\infty. 
\]
\end{lem}
\begin{proof}
Let $O$ be another open ball such that $\overline{U}\subset O$. Define $h_n>0$ by 
\[ 
 h_n^{-q}= \max\left\{n, \max_{1\le k\le n}\|v_k\|_{C^{2,3}([0,T)\times O)}\right\}. 
\]
Then $h_n$ is decreasing and satisfies $\lim_{n\to\infty}h_n=0$. Clearly we have 
\[
 \|v_n\|_{C^{2,3}([0,T)\times O)}\le h_n^{-q}
\]
as well as 
\[
 |D_xv_n(t,x)| + |D_x^2v_n(t,x)|\le c_1 h_n^{-q}, \quad 0\le t<T, \;\; x\in O, 
\]
for some $c_1>0$. 

Let $\{t_i\}_{i=1}^N$ be such that $0=t_0<t_1<\cdots <t_N=T$ and $\max_i(t_i -t_{i-1})\le h_n$. Then the condition (A4) means 
\[
|v_n(t_N,x)|\le B_n, \quad x\in U, \;\; n\ge 1 
\]
for some positive constant $B_n$.  
So suppose that for $i\le N-1$ there exists $B_{i+1}>0$ such that
\[
 |v_n(t_{i+1},x)|\le B_{i+1}, \quad x\in U, \;\; n\ge n_1
\]
with some $n_1\ge 1$ to be determined below. 
Note that we have 
\[
 \|v_n(t_i)\|_{C^3(O)}\le h_n^{-q}, \quad i=0,1,\ldots,N-1. 
\] 
By a routine argument, we can extend each $v_n(t_i,\cdot)$ to a $C^3$-function on $\mathbb{R}^d$. Indeed, 
recall that there exists a $C^{\infty}$-function $\zeta$ compactly supported in $O$ such that 
$0\le\zeta\le 1$ on $\mathbb{R}^d$ and $\zeta=1$ on $U$. See, e.g., 
Theorem 1.4.1 in H{\"o}rmander \cite{hor:1990}. 
Then, for $i=0,1,\ldots,N-1$ define $u_n(t_i,\cdot)\in C^3(\mathbb{R}^d)$ by 
\[
 u_n(t_i,x) = v_n(t_i,x)\zeta(x) + h_n^{-q}(1-\zeta(x)), \quad x\in\mathbb{R}^d. 
\]
It is straightforward to see that  
$D^{\alpha}_xu_n(t_i,x)=D^{\alpha}_xv_n(t_i,x)$ for $x\in U$, $\|u_n(t_i)\|_{C^3(\mathbb{R}^d)}\le c_2h_n^{-q}$ and 
$|D_xu_n(t_i,x)| + |D_x^2u_n(t_i,x)|\le c_2h_n^{-q}$, $x\in \mathbb{R}^d$, 
for some positive constant $c_2$. 

To get a bound of $v_n(t_i,\cdot)$, put 
\[
 R_n(x)=\frac{1}{h_n}(v_n(t_i,x) - v_n(t_{i+1},x)) + F(t_{i+1},x; v_n(t_{i+1})).   
\]
Then, 
\begin{align*}
|R_n(x)| &= \left|-\frac{1}{h_n}\int_{t_i}^{t_{i+1}}\partial_tv_{n}(t,x)\, dt + F(t_{i+1}, x; v_{n}(t_{i+1},\cdot))\right| \\ 
&\le |-\partial_tv_{n}(t_i,x)+F(t_{i+1},x; v_{n}(t_{i+1},\cdot))| 
   + \frac{1}{h_n}\int_{t_i}^{t_{i+1}}|\partial_tv_{n}(t,x) - \partial_tv_{n}(t_i,x)|\,dt. 
\end{align*}
By (B1), the first term of the right-hand side in the inequality just above converges to zero as $n\to\infty$. 
By the definition of $h_n$, the second term is at most $Ch_n^{1-q}$. Hence $R_n(x)$ is bounded with respect to $n$.  
Thus Lemma \ref{lem:2.2} yields  
$|v_n(t_i,x)|\le |W| + Ch$ where $x\in U$ and 
\begin{equation*}
W=\sup_{(p,\Gamma)\in\mathcal{D}_{h,c_2}}\inf_{w\in\mathcal{X}_{h,\kappa}}
  \left[u_n(t_{i+1},x+\sqrt{h}w) 
 -\sqrt{h}w^{\mathsf{T}}p-\frac{h}{2}w^{\mathsf{T}}\Gamma w  
 -hF(t,x,v_n(t_{i+1},x),p,\Gamma)\right]
\end{equation*}
with $h=h_n$. 
Considering $p=0$ and $\Gamma=0$, we see 
$W\ge -(1+C_0h)B_{i+1}-C_0h$. 

To obtain an upper bound, observe 
\[
 W\le B_{i+1}+  
  \sup_{(p,\Gamma)\in\mathcal{D}_{h,c_2}}\inf_{w\in\mathcal{X}_{h,\kappa}}
   W_{p,\Gamma ,w}, 
\]
where 
\[
W_{p,\Gamma,w}= -\sqrt{h}w^{\mathsf{T}}p-\frac{h}{2}w^{\mathsf{T}}\Gamma w 
    -hF(t,x,v_n(t_{i+1},x),p,\Gamma).
\] 
Then we will show that for any $(p,\Gamma)\in\mathcal{D}_{h,c_2}$ 
we can find $w\in\mathcal{X}_{h,\kappa}$ satisfying 
$W_{p,\Gamma,w}\le K_1hB_{i+1} + Ch$. 
So fix $(p,\Gamma)\in\mathcal{D}_{h,c_2}$. 
First assume that the minimum eigenvalue of $-\Gamma$ is 
greater than or equal to $-h^{\gamma}$, where $\gamma\in (0,2\kappa-q)$. 
If $p=0$ then we may take $w=0$, leading to 
$W_{p,\Gamma,w}\le -hF(t,x,v_n(t_{i+1},x),0,h^{\gamma}I) 
\le C_0h+C_0hB_{i+1}+C_0h^{1+\gamma}$. 
Otherwise, take $w=h^{q}p/|p|$. Then we see 
\begin{align*}
W_{p,\Gamma,w}&\le -h^{(1/2)+q}|p|+\frac{h^{1+2q+\gamma}}{2} + 
 C_0h(1+B_{i+1}) + C_0h|p| + C_0h^{1+\gamma} \\
 &\le |p|(-h^{(1/2)+q} + C_0 h) +C_0h+ C_0hB_{i+1} \le Ch+C_0hB_{i+1}
\end{align*}
since $-h^{2/3}+C_0h\le 0$ for $n\ge n_1$ 
with some $n_1\in\mathbb{N}$.  

Next assume that the minimum eigenvalue of $-\Gamma$ is less than 
$-h^{\gamma}$. Then take $w$ to be the corresponding eigenvector satisfying 
$-p^{\mathsf{T}}w\le 0$ and $|w|=h^{-\kappa}$. 
This choice leads to 
\begin{equation*}
 W_{p,\Gamma,w}\le -\frac{h^{1+\gamma-2\kappa}}{2}+ C_0h(1+B_{i+1}) 
  +2ch^{1-q} \le C_0h(1+B_{i+1})
\end{equation*}
since there exists $h_1\in (0,h_1^{\prime}]$ such that 
$-h^{1+\gamma-2\kappa+q} + 4c h\le 0$ for 
$n\ge 1$. 

Therefore we deduce that $|W|\le (1+C_0h)B_{i+1}+Ch$ for $n\ge n_1$. 
Denoting the right-hand side by $B_i$, we obtain the sequence $\{B_i\}$ 
satisfying $B_i=(1+C_0h)B_{i+1}+Ch$. By a routine argument we have 
$B_i\le e^{TC_0}B_n + Ce^{TC_0}$ for all $i$. 

Hence we have $\max_{i=0,1,\ldots,N}\sup_{x\in U}|v_n(t_i,x)|\le C$. 
By (B1), for any $t$ we can take a nearest $t_i$ to get $|v_n(t,x)|\le |v_n(t_i,x)| + Ch_n^{1-q}$. 
Thus the lemma follows. 
\end{proof}

\begin{proof}[Proof of Theorem $\ref{thm:2.1}$]
The argument of the proof is similar to that in \cite{nak:2017}. We will show that  
\begin{equation*}
 \overline{v}(t,x)=\limsup_{{s\to t, \; y\to x}
  \atop{n\to\infty}} v_n(s,y), 
 \quad (t,x)\in [0,T]\times\mathbb{R}^d, 
\end{equation*}
is a viscosity subsolution of (\ref{eq:1.1}). 
Notice that $\overline{v}$ is finite on $[0,T]\times\mathbb{R}^d$ by Lemma \ref{lem:2.3}. 

Fix $(t,x)\in [0,T)\times\mathbb{R}^d$ and 
let $\varphi$ be a $C^3$-function on $[0,T]\times\mathbb{R}^d$ such that 
$\overline{v}-\varphi$ has a global strict maximum at $(t,x)$ with $(\overline{v}-\varphi)(t,x)=0$. 
By definition of $\overline{v}$, there exist $n_m$, $t_m$, $y_m$ such that 
as $m\to \infty$, 
\[
 n_m\to \infty, \;\; (t_{m},y_m)\to (t,x), \;\; 
 v_{n_m}(t_{m},y_m)\to \overline{v}(t,x). 
\]
and that 
\begin{equation}
\label{eq:2.2}
  c_m:=(v_{n_m}-\varphi)(t_{m},y_m)\ge \sup_{(s,y)\in [0,T)\times\mathbb{R}^d}
 (v_{n_m}-\varphi)(s,y) - h_m^{3/2}. 
\end{equation}
In particular, $c_m\to 0$. 
Take an open neighborhood $O$ of $x$. 
It follows from (\ref{eq:2.2}) that for any $y\in O$ we have 
\begin{equation}
\label{eq:2.3}
 \varphi(t_{m}+h_m,y)+c_m + h_m^{3/2}
 \ge v_{n_m}(t_{m}+h_m,y), 
\end{equation}
where 
\[ 
 h_m^{-q}= \max\left\{n_m, \max_{\nu\le m}\|v_{n_{\nu}}\|_{C^{2,3}([0,T)\times O)}\right\}
\]
and $q$ is as in Lemma \ref{lem:2.2}. 

Now rewrite $v_{n_m}(t_{m},y_m)$ as
\begin{equation}
\label{eq:2.4p}
 v_{n_m}(t_{m},y_m)
  = v_{n_m}(t_{m} + h_m,y_m)
  - h_m F(t_{m}+h_m,y_m;v_{n_m}(t_{m}+h_m,\cdot)) + h_m J_m. 
\end{equation}
Then, 
\begin{align*}
|J_m| &= \left|-\frac{1}{h_m}\int_{t_m}^{t_m+h_m}\partial_tv_{n_m}(t,y_m)\, dt + F(t_m, y_m; v_{n_m}(t_m,\cdot))\right| \\ 
&\le |-\partial_tv_{n_m}(t_m,y_m)+F(t_m,y_m; v_{n_m}(t_m,\cdot))| 
   + \frac{1}{h_m}\int_{t_m}^{t_m+h_m}|\partial_tv_{n_m}(t,y_m) - \partial_tv_{n_m}(t_m,y_m)|\,dt. 
\end{align*}
By (B1), the first term of the right-hand side in the inequality just above converges to zero as $m\to\infty$. 
By the definition of $h_m$, the second term is at most $Ch_m^{1-q}$. Hence we have $\lim_{m\to\infty}J_m=0$. 
 
Again by the definition of $h_m$ there exists a constant $c_3>0$ such that 
$\sum_{|\alpha|_1\le 3}\sup_{x\in O}|D^{\alpha}\phi(x)|\le c_3 h^{-q}$ for $\phi=v_{n_m}(t_{m}), \varphi$, $m\ge 1$. 
With the representation (\ref{eq:2.4p}), we apply Lemma \ref{lem:2.2} and use the inequality (\ref{eq:2.3}) to get that 
there exist $\kappa, \beta>0$ and $m_0\in\mathbb{N}$ such that for any $m\ge m_0$,  
\begin{align*}
 &v_{n_m}(t_{m},y_m) \\
 &\le 
  \sup_{(p,\Gamma)\in\mathcal{D}_{h_m, c_4}}
  \inf_{w\in\mathcal{X}_{h_m,\kappa}}\Big[
 v_{n_m}(t_{m}+h_m,y_m+\sqrt{h_m}w)-\sqrt{h_m}p^{\mathsf{T}}w 
 -\frac{h_m}{2}w^{\mathsf{T}}\Gamma w \\
 &\qquad - h_mF(t_{m}+h_m,y_m,v_{n_m}(t_{m}+h_m,y_m),p,\Gamma)\Big] 
  + h_m J_m + Ch_m^{1+\beta} \\
 &\le \sup_{p,\Gamma}\inf_w\Big[
  \varphi(t_{m}+h_m,y_m+\sqrt{h_m}w)-\sqrt{h_m}p^{\mathsf{T}}w 
 -\frac{h_m}{2}w^{\mathsf{T}}\Gamma w \\
 &\qquad -h_mF(t_{m}+h_m,y_m,v_{n_m}(t_{m}+h_m,y_m),p,\Gamma)\Big] 
  +c_m + h_m^{3/2}+h_m J_m + Ch_m^{1+\beta} \\
 &\le \varphi(t_{m}+h_m,y_m) 
         -h_m F(t_{m}+h_m,y_m,v_{n_m}(t_{m}+h_m,y_m),
     D\varphi(t_{m}+h_m,y_m),D^2\varphi(t_{m}+h_m,y_m)) \\
 &\qquad +c_m+h_m^{3/2}
    +h_m J_m  + Ch_m^{1+\beta}
\end{align*}
for a suitable constant $c_4>0$. 
This together with $v_{n_m}(t_{m},y_m)=c_m+\varphi(t_{m},y_m)$ leads to 
\begin{equation*}
\begin{aligned}
 &-\frac{1}{h_m}\left(\varphi(t_{m}+h_m,y_m)-\varphi(t_{m},y_m)\right) \\
 &\quad +F(t_{m},y_m,v_{n_m}(t_{m},y_m),D\varphi(t_{m},y_m),
    D^2\varphi(t_{m},y_m))
 \le o(1)
\end{aligned} 
\end{equation*}
for any sufficiently large $m$.  
Sending $m\to\infty$, we have 
\begin{equation*}
 -\partial_t\varphi(t,x)+F(t,x,\overline{v}(t,x),D\varphi(t,x),D^2\varphi(t,x)) 
 \le 0, 
\end{equation*}
whence the subsolution property at $(t,x)$. 

In the case $(t,x)\in\{T\}\times\mathbb{R}^d$, from (B2) we have 
$\overline{v}(t,x)=g(x)$. Therefore $\overline{v}$ is a viscosity subsolution of \eqref{eq:1.1}.

A similar argument shows that 
\[
\underline{v}(t,x)=\liminf_{{s\to t,\; y\to x}
 \atop{n\to\infty}}v_n(s,y), 
\quad (t,x)\in [0,T]\times\mathbb{R}^d
\]
is a viscosity supersolution of (\ref{eq:1.1}). 
By (A5), we obtain $\overline{v}\le \underline{v}$. 
This and $\overline{v}\ge \underline{v}$ means $\overline{v}=\underline{v}$. 
From this the conclusion of the theorem follows.  
\end{proof}

\subsection{Hamilton--Jacobi--Bellman equations}\label{sec:2.2}

Here we consider the case of HJB equations. Precisely, $F$ is assumed to be given by 
\begin{equation}
\label{eq:2.4}
 F(t,x,z,p,X)= \sup_{a\in A}\left[-b(t,x,a)^{\mathsf{T}}p - \frac{1}{2}\mathrm{tr}((\sigma\sigma^{\mathsf{T}})(t,x,a)X) - f(t,x,a)\right], 
\end{equation}
for $t\in [0,T]$, $x\in\mathbb{R}^d$, $z\in\mathbb{R}$, $p\in\mathbb{R}^d$, and $X\in\mathbb{S}^d$, where $A\subset \mathbb{R}^{d_1}$ is compact,
with some continuous functions $b:[0,T]\times\mathbb{R}^d\times A\to\mathbb{R}^d$,
$\sigma: [0,T]\times\mathbb{R}^d\times A\to\mathbb{R}^{d\times d_2}$, and $f:[0,T]\times\mathbb{R}^d\times A\to\mathbb{R}$.

Thus in this case \eqref{eq:1.1} becomes the HJB equation
\[
\left\{
\begin{split}
 & \partial_t v + \inf_{a\in A}\mathcal{L}^av(t,x)=0, \quad (t,x)\in 
 [0,T)\times \mathbb{R}^d, \\
 &v(T,x)=g(x), \quad x\in\mathbb{R}^d, 
\end{split}
\right. 
\]
where for $a\in A$,  
\[
 \mathcal{L}^a\phi(t,x)=b(t,x,a)^{\mathsf{T}}D_x\phi(t,x) + \frac{1}{2}\mathrm{tr}(\sigma\sigma^{\mathsf{T}}(t,x,a)D^2_{xx}\phi(t,x)) + f(t,x,a), \quad (t,x)\in [0,T)\times\mathbb{R}^d. 
\] 

Suppose that as in Section \ref{sec:2.1} a candidate of approximate solution $v_n$ satisfies $\|v_n\|_{C^{2,3}([0,T]\times\mathbb{R}^d)}<\infty$ for any $n\ge 1$. 
Then choose a positive sequence $\{\lambda_n\}_{n=1}^{\infty}$ satisfying 
\[
 \|v_n\|_{C^{1,2}([0,T]\times\mathbb{R}^d)}\le \lambda_n, \quad n\ge 1. 
\]
Let $O_n\subset \mathbb{R}^d$ be open and contain an open ball with radius $R_n>0$ centered at a given $\bar{x}\in\mathbb{R}^d$, i.e., 
$O_n\supset\{x\in\mathbb{R}^d: |x-\bar{x}|<R_n\}$. We assume that $\lim_{n\to\infty}R_n=\infty$. Put 
\[
 \eta_n=\sup_{(t,x)\in [0,T]\times O_n}\left|\inf_{a\in A}(\partial_t+\mathcal{L}^{a})v_n(t,x)\right| \vee 
   \sup_{x\in O_n}|v_n(T,x) - g(x)|
\]
for any $n\ge 1$, where $a\vee b=\max\{a,b\}$. 

To derive the error bound between $v_n$ and the true solution $v$, we consider the following condition: 
\begin{enumerate}
\item[(A6)] There exists a constant $C_1>0$ such that for $\psi=b,\sigma,f$, and for $t,t^{\prime}\in [0,T]$, $x,x^{\prime}\in\mathbb{R}^{d}$, $a\in A$, 
  \begin{align*}
   |\psi(\xi,a)|&\le C_1, \\  
   |\psi(\xi,a)-\psi(\xi^{\prime},a)|&\le C_1|\xi-\xi^{\prime}|. 
  \end{align*}
\end{enumerate}
It is well-known that under (A4) and (A6) there exists a unique continuous viscosity solution $v$ of \eqref{eq:1.1} (see, e.g., Fleming and Soner \cite{fle-son:2006}). 

\begin{thm}
\label{thm:2.2}
Suppose that \rm{(A4)} and \rm{(A6)} hold.
Suppose moreover that $\|v_n\|_{C^{2,3}([0,T]\times\mathbb{R}^d)}<\infty$ and the functions
$\partial_tv_n$, $D_xv_n$, $D_x^2v_n$ are all Lipschitz on $[0,T]\times\mathbb{R}^d$ for each $n\ge 1$.
Then for any compact set $Q\subset [0,T]\times\mathbb{R}^d$ there exist a constant $C>0$ and $n_0\in\mathbb{N}$ such that for $n\ge n_0$ 
\[
 \sup_{(t,x)\in Q}|v(t,x)-v_n(t,x)|\le (1+T)\eta_n + C \lambda_n\exp\left(- \frac{R_n^2}{8C_1^2T}\right). 
\]
\end{thm}
\begin{proof}
It is known that $v$ has the following stochastic control representation 
\begin{equation}
\label{eq:2.5}
 v(t,x)=\inf_{\alpha\in\mathcal{A}}\mathbb{E}\left[\int_t^Tf(s,X_s^{t,x,\alpha},\alpha_s)\,ds + g(X_T^{t,x,\alpha})\right], 
 \quad (t,x)\in[0,T]\times\mathbb{R}^d. 
\end{equation}
See, e.g., Pham \cite[Chapter 4]{pha:2009}. Here, $\{X_s^{t,x,\alpha}\}_{t\le s\le T}$ is a unique solution of the controlled stochastic differential equation 
\[
 dX_s^{t,x,\alpha}=b(s,X_s^{t,x,\alpha},\alpha_s)\,ds +\sigma(s,X_s^{t,x,\alpha},\alpha_s)\, dW_s 
\]
with initial condition $X_t^{t,x,\alpha}=x$ on a complete probability space $(\Omega,\mathcal{F},\mathbb{P})$ with a filtration $\mathbb{F}=\{\mathcal{F}_s\}_{s\ge 0}$ 
satisfying the usual conditions.  
The process $\{W_s\}_{s\ge 0}$ is a $d_2$-dimensional standard $\mathbb{F}$-Brownian motion.
The process $\{\alpha_s\}_{0\le s\le T}$ is an $A$-valued $\mathbb{F}$-progressively measurable process. 
We have denoted by $\mathcal{A}$ the collection of such $\alpha$'s. 

Since $\|v_n\|_{C^{2,3}([0,T]\times\mathbb{R}^d)}<\infty$ for any $n\ge 1$, the function 
\[
 r_n(t,x):= - \partial_tv_n(t,x) - \inf_{a\in A}\mathcal{L}^av_n(t,x), \quad (t,x)\in [0,T]\times\mathbb{R}^d, 
\]
is bounded and Lipschitz on $[0,T]\times\mathbb{R}^d$. From this and the trivial equality 
\[
 \partial_tv_n(t,x) + \inf_{a\in A}\mathcal{L}^{a}v_n(t,x) + r_n(t,x)=0
\]
we again use the stochastic control representation to get 
\[
 v_n(t,x) = \inf_{\alpha\in\mathcal{A}}\mathbb{E}\left[\int_t^T(f(s,X_s^{t,x,\alpha},\alpha_s) + r_n(s,X_s^{t,x,\alpha}))\,ds + v_n(T, X_T^{t,x,\alpha})\right], 
 \quad (t,x)\in [0,T]\times\mathbb{R}^d. 
\]
This together with \eqref{eq:2.5} leads to 
\begin{align*}
&|v(t,x)-v_n(t,x)| \\ 
&\le \sup_{\alpha\in\mathcal{A}}\int_t^T\mathbb{E}|r_n(s,X_s^{t,x,\alpha})|\, ds + \sup_{\alpha\in\mathcal{A}}\mathbb{E}|g(X_T^{t,x,\alpha})-v_n(T,X_T^{t,x,\alpha})| \\
&\le \sup_{\alpha\in\mathcal{A}}\int_t^T\mathbb{E}|r_n(s,X_s^{t,x,\alpha})|1_{\{X_s^{t,x,\alpha}\in O_n\}}\, ds
 + \sup_{\alpha\in\mathcal{A}}\int_t^T\mathbb{E}|r_n(s,X_s^{t,x,\alpha})|1_{\{X_s^{t,x,\alpha}\notin O_n\}}\, ds \\
&\quad + \sup_{\alpha\in\mathcal{A}}\mathbb{E}|g(X_T^{t,x,\alpha})-v_n(T,X_T^{t,x,\alpha})|1_{\{X_s^{t,x,\alpha}\in O_n\}} 
             + \sup_{\alpha\in\mathcal{A}}\mathbb{E}|g(X_T^{t,x,\alpha})-v_n(T,X_T^{t,x,\alpha})|1_{\{X_s^{t,x,\alpha}\notin O_n\}} \\ 
&\le (1+T)\eta_n + C\|v_n\|_{C^{1,2}([0,T]\times\mathbb{R}^d)} \sup_{\alpha\in\mathcal{A}}\mathbb{P}\left(\sup_{t\le s\le T}|X_s^{t,x,\alpha}-x_0|\ge R_n\right). 
\end{align*}

Now we will show that 
\begin{equation}
\label{eq:2.6}
 \sup_{(t,x)\in Q}\sup_{\alpha\in\mathcal{A}}\mathbb{P}\left(\sup_{t\le s\le T}|X_s^{t,x,\alpha}-x_0|\ge R_n\right) 
 \le \exp\left(- \frac{R_n^2}{8C_1^2T}\right), \quad n\ge n_0
\end{equation}
for some $n_0\ge 1$, 
which completes the proof of the theorem. 

Fix $t,x,\alpha$ and put $b_s=b(s,X_s^{t,x,\alpha},\alpha_s)$ and $\sigma_s=\sigma(s,X_s^{t,x,\alpha},\alpha_s)$. Observe 
\[
 \left\{\sup_{t\le s\le T}|X_s^{t,x,\alpha} - x_0|\ge R_n\right\} \subset \left\{|x_0| + |x| + \int_t^T|b_s|\,ds\ge \frac{R_n}{2}\right\} 
 \bigcup \left\{\sup_{t\le s\le T}\left|\int_t^T\sigma_s\,dW_s\right| \ge \frac{R_n}{2}\right\}. 
\]
Since $R_n\to \infty$, there exists $n_0\ge 1$ such that $R_n>2(|x_0| + \sup_{(s,y)\in Q}|y| + C_1T)$ for $n\ge n_0$. Thus
\[
 \mathbb{P}\left(\sup_{t\le s\le T}|X_s^{t,x,\alpha}-x_0|\ge R_n\right)\le \mathbb{P}\left(\sup_{t\le s\le T}\left|\int_t^T\sigma_s\,dW_s\right|\ge \frac{R_n}{2}\right).
\]
Consider the stopping time $\tau=\inf\{s\in [t,T]: |\int_t^s\sigma_u\,dW_u|\ge R_n/2\}\wedge T$. 
Since $\sigma$ is bounded and $\int_t^s\sigma_u\,dW_u$ is a continuous martingale, the optional sampling theorem means 
\begin{align*}
 1&=\inf_{\mu>0}\mathbb{E}\left[\exp\left(\mu\int_t^{\tau}\sigma_u\,dW_u - \frac{\mu^2}{2}\int_t^{\tau}|\sigma_u|^2\,du\right)\right] \\ 
  &\ge \inf_{\mu>0}\mathbb{E}\left[1_{\{\sup_{t\le s\le T}|\int_t^s\sigma_u\,dW_u|\ge R_n/2\}} 
        \exp\left(\mu\int_t^{\tau}\sigma_u\,dW_u - \frac{\mu^2}{2}\int_t^{\tau}|\sigma_u|^2\,du\right) \right] \\ 
  &\ge \mathbb{P}\left(\sup_{t\le s\le T}\left|\int_t^s\sigma_u\,dW_u\right|\ge \frac{R_n}{2}\right)
    \exp\left(\min_{\mu>0}\left(\frac{R_n\mu}{2} - \frac{C_1^2T\mu^2}{2}\right)\right) \\ 
  &= \mathbb{P}\left(\sup_{t\le s\le T}\left|\int_t^s\sigma_u\,dW_u\right|\ge \frac{R_n}{2}\right)
        \exp\left(\frac{R_n^2}{8C_1^2T}\right)
\end{align*}
for $n\ge n_0$. Thus \eqref{eq:2.6} follows. 
\end{proof}

\section{Applications to kernel-based function approximation}\label{sec:3}

Here we shall apply our abstract result in Section \ref{sec:2} to some function approximation methods. 
We seek a solution of \eqref{eq:1.1} that is represented as a linear combination of given basis functions. 

\subsection{Convergence for general equations}\label{sec:3.1}

This subsection presents kernel-based function approximation
methods that illustrate how the abstract convergence framework developed in
Section \ref{sec:2.1} applies to differentiable and non-monotone representations.

The formulation is designed to provide a concrete setting in which the consistency
and stability assumptions \textnormal{(B1)--(B2)} can be verified within a
kernel-based framework.

Let $\Phi: \mathbb{R}^{d+1}\to \mathbb{R}$ be a radial and positive definite function, 
i.e., $\Phi(\cdot)=\phi(|\cdot|)$ for some $\phi:[0,\infty)\to\mathbb{R}$ and 
for every $\ell\in\mathbb{N}$, for all pairwise distinct 
$y_1,\ldots, y_{\ell}\in\mathbb{R}^{d+1}$ and for all 
$\alpha=(\alpha_i)\in\mathbb{R}^{\ell}\setminus\{0\}$, we have 
\begin{equation*}
  \sum_{i,j=1}^{\ell}\alpha_i\alpha_j\Phi(y_i-y_j)>0. 
\end{equation*}
Let the Fourier transform $\widehat{\Phi}$ of $\Phi$ satisfy  
\[
 c_1(1+|\xi|^2)^{-\tau}\le \widehat{\Phi}(\xi) \le c_2(1+|\xi|^2)^{-\tau}, \quad \xi\in\mathbb{R}^{d+1}, 
\]
where $c_1,c_2$ are positive constants with $c_1\le c_2$, and $\tau\in\mathbb{N}$ with $\tau>3 + (d+1)/2$. 
Then, there exists a unique Hilbert space $\mathcal{H}\subset C^3((0,T)\times\mathbb{R}^{d})$ with 
norm $\|\cdot\|_{\mathcal{H}}$, called the 
native space, of real-valued functions on $(0,T)\times \mathbb{R}^{d}$ such that 
$\Phi$ is a reproducing kernel for $\mathcal{H}$. 
Moreover, 
$\mathcal{H}$ is $L^2$-Sobolev space on $(0,T)\times\mathbb{R}^{d}$ of order $\tau$ with equivalent norm, and 
$\mathrm{span}\{\Phi(\cdot - \xi); \xi\in (0,T)\times \mathbb{R}^{d}\}$ is dense in $\mathcal{H}$. 
See, e.g., Wendland \cite[Chapter 10]{wen:2010}.

Let $Q_n=(0,T)\times \{x\in\mathbb{R}^d: |x-\bar{x}|_{\infty}<R_n\}$ where $\bar{x}\in\mathbb{R}^{d}$, $R_n>0$. 
Let $Q=\{y\in\mathbb{R}^{d+1}: |y|_{\infty}<1\}$.  
Take a pairwise distinct points set $\Gamma:=\{\tilde{y}^1,\ldots,\tilde{y}^N\}\subset Q$, where $N=N_n\in\mathbb{N}$ also depends on $n$. Then 
define $y^j=\iota(\tilde{y}^j)$ and $\Gamma_n=\{y^1,\ldots,y^N\}\subset Q_n$ where 
\[
 \iota(y)=\left(\frac{T}{2}(y_0+1), R(y_1+\bar{x}_1),\ldots, R(y_d+\bar{x}_d)\right)^{\mathsf{T}}\in Q_n, \quad y=(y_0,y_1,\ldots,y_d)\in Q. 
\] 
Further, put $t^j=y^j_0$ and $x^j=(y_1^j,\ldots,y_d^j)$. 
Then we shall construct an approximate solution $v_n$ of the form 
\begin{equation}
\label{eq:3.1}
 v_n(y) = \sum_{j=1}^N\theta_j\Phi(y-y^j), \quad y=(t,x)\in [0,T]\times \mathbb{R}^d. 
\end{equation}
Substituting this into \eqref{eq:1.1} and evaluating the resulting equality at a given $\Gamma_n^e=\{\tilde{y}^1, \ldots, \tilde{y}^M\}\subset Q_n$, 
we are led to the following finite-dimensional optimization problem: to minimize
\begin{equation}
\label{eq:3.2}
 \max_{\ell=1,\ldots,M}\left|-\sum_{j=1}^N\theta_j\partial_t\Phi(\tilde{y}^{\ell}-y^j) + F\left(\tilde{y}^{\ell}; \sum_{j=1}^N\theta_j\Phi((\tilde{t}^{\ell},\cdot) - y^j)\right)\right| 
 \vee \left|\sum_{j=1}^N\theta_j\Phi((T,\tilde{x}^{\ell})-y^j) - g(\tilde{x}^{\ell})\right| 
\end{equation}
over all $\theta\in\mathbb{R}^N$, where $a\vee b=\max(a,b)$ for $a,b\in\mathbb{R}$ and $(\tilde{t}^{\ell},\tilde{x}^{\ell})=\tilde{y}^{\ell}$. 
We assume $\Gamma^e_n\supset\Gamma_n$. 

To meet (B1), we assume that there exists a classical solution $v$ of \eqref{eq:1.1} such that $v\in\mathcal{H}$.
Thus we impose
\[
 \|v_n\|_{\mathcal{H}}^2= \theta^{\mathsf{T}}K\theta\le \lambda^2,
\]
where $K=\{\Phi(\tilde{y}^{\ell} - y^j)\}_{{1\le\ell\le M}\atop{1\le j\le N}}\in\mathbb{R}^{M\times N}$.

We shall adopt an epigraph formulation of the problem \eqref{eq:3.2}. 
Consequently, the approximate solution is constructed by solving the following problem: 
\begin{equation}
\tag{$P_{\lambda}$}
\left\{
\begin{aligned}
 &\text{minimize}\;\; \gamma, \\
 &\text{s.t.}\;\;  (\theta, \gamma) \in \mathbb{R}^N\times [0,\infty), \;\; \theta^{\mathsf{T}}K\theta\le \lambda^2, \;\;
  -\gamma \mathbf{1}_M\le \tilde{K}\theta - \tilde{g} \le \gamma \mathbf{1}_M, \;\;
  -\gamma \mathbf{1}_M\le \mathbf{F}(\theta) \le \gamma \mathbf{1}_M,
\end{aligned}
\right.
\end{equation}
where $\tilde{K}=\{\Phi((T,\tilde{x}^{\ell}) - y^j)\}_{{1\le\ell\le M}\atop{1\le j\le N}}$, $\tilde{g}=(g(\tilde{x}^1),\ldots,g(\tilde{x}^N))^{\mathsf{T}}$, and 
$\mathbf{F}(\theta)=(\mathbf{F}^1(\theta),\ldots,\mathbf{F}^M(\theta))^{\mathsf{T}}$ with 
\[
 \mathbf{F}^{\ell}(\theta)= - \sum_{j=1}^N\theta_j(\partial_t\Phi)(\tilde{y}^{\ell}-y^j) + F\left(\tilde{y}^{\ell}; \sum_{j=1}^N\theta_j\Phi((\tilde{t}^{\ell},\cdot) - y^j)\right), 
 \quad \ell=1,\ldots,M, \;\; \theta\in\mathbb{R}^N. 
\]
The inequalities and function evaluations in ($P_{\lambda}$) should be understood in the pointwise sense. 

Our candidate of approximate solution is now constructed as follows: 
given $\lambda$ and $\varepsilon=\varepsilon_n>0$, take an $\varepsilon_n$-optimal solution $(\theta^{(n)},\gamma^{(n)})\in\mathbb{R}^N\times [0,\infty)$
of ($P_{\lambda}$). Then, define $v_n^{(\lambda)}\in\mathcal{H}$ by 
\eqref{eq:3.1} with $\theta=\theta^{(n)}$.

We shall give a convergence guarantee of the methods above using the result in Section \ref{sec:2}. 
To this end, we consider the following condition: 
\begin{enumerate}
\item[(A1$^{\prime}$)] There exists a constant $C_2>0$ such that for $t,t^{\prime}\in [0,T]$, $x,x^{\prime}\in\mathbb{R}^d$, $z,z^{\prime}\in\mathbb{R}$,
$p,p^{\prime}\in\mathbb{R}^d$, $X,X^{\prime}\in\mathbb{S}^d$, 
\begin{align*}
 &|F(t,x,z,p,X) - F(t^{\prime},x^{\prime},z^{\prime},p^{\prime},X^{\prime})| \\
 &\le C_2(1+|p|+|p^{\prime}| + |X| + |X^{\prime}|)(|t-t^{\prime}| + |x-x^{\prime}| + |z-z^{\prime}|) +C_2(|p-p^{\prime}| + |X-X^{\prime}|). 
\end{align*}
\end{enumerate}
Let $h_n$ be the fill distance between $Q$ and $\Gamma$, defined by 
\[
 h_n=\sup_{y\in Q}\min_{j=1,\ldots,N_n}|y-y^j|. 
\]
It is straightforward to see that the fill distance between $Q_n$ and $\Gamma_n$ is $O(R_nh_n)$. 
Then we have the following: 
\begin{thm}
\label{thm:3.1}
Suppose that \rm{(A1$^{\prime}$)}, \rm{(A2)}--\rm{(A5)} hold and that 
there exists a classical solution $v\in C([0,T]\times\mathbb{R}^d)\cap\mathcal{H}$ of $\eqref{eq:1.1}$. 
Suppose moreover that 
\[
 \varepsilon_n\to 0, \quad R_nh_n\to 0, \quad n\to\infty.  
\]
Then, there exists $\lambda>0$ such that $\{v_n^{(\lambda)}\}$ satisfies the consistencies \rm{(B1)} and \rm{(B2)}. In particular, 
\[
v_n^{(\lambda)}(s,x)\to v(t,x),   
\]
as $s\to t$ and $n\to\infty$ uniformly on any compact subset of $\mathbb{R}^d$. 
\end{thm}
\begin{proof}
Step (i). 
Let $\lambda>0$ be fixed, to be determined below. In view of the Sobolev inequality and the constraints in the problem (P$_{\lambda}$), we find that 
\begin{equation}
\label{eq:3.3}
 \|v_n^{(\lambda)}\|_{C^{3}([0,T]\times \mathbb{R}^{d})}\le C\lambda.
\end{equation} 
Fix $(t,x)\in (0,T)\times\mathbb{R}^d$ and let $(t_n,x_n)\in (0,T)\times\mathbb{R}^d$, $n\ge 1$, be an arbitrary sequence such that $y_n=(t_n,x_n)\to (t,x)$ 
as $n\to\infty$. Since $R_n\to\infty$ there exists $n_0\ge 1$ such that $y=(t,x)\in Q_n$ for any $n\ge n_0$. 
Then take a nearest $y^*=(t^*,x^*)\in\Gamma_n$ in the sense that 
$|y-y^*|=\min_{j=1,\ldots,N}|y-y^j|$. The condition \rm{(A1$^{\prime}$)} and \eqref{eq:3.3} yield
\begin{align*}
&|-\partial_tv_n^{(\lambda)}(t_n,x_n) + F(t_n,x_n; v_n^{(\lambda)}(t_n,\cdot)) + \partial_tv_n^{(\lambda)}(t,x) - F(t,x; v_n^{(\lambda)}(t,\cdot))| \\
&\le |\partial_tv_n^{(\lambda)}(y_n) - \partial_tv_n^{(\lambda)}(y)| 
 + C_2(1+|Dv_n^{(\lambda)}(y)| + |Dv_n^{(\lambda)}(y_n)| + |D^2v_n^{(\lambda)}(y)| + |D^2v_n^{(\lambda)}(y_n)|) \\ 
 &\hspace*{13em} \times (|y-y_n| + |v_n^{(\lambda)}(y_n)-v_n^{(\lambda)}(y)|) \\ 
 &\quad + C_2(|Dv_n^{(\lambda)}(y_n)-Dv_n^{(\lambda)}(y)| + |D^2v_n^{(\lambda)}(y_n) - D^2v_n^{(\lambda)}(y)|) \\ 
 &\le C(1+\lambda+ \lambda^2)(|y-y_n|) \\ 
 &\le C(1+\lambda+ \lambda^2)R_nh_n. 
\end{align*}
Similarly, 
\[
 |-\partial_tv_n^{(\lambda)}(y) + F(y_n; v_n^{(\lambda)}(t_n,\cdot)) + \partial_tv_n^{(\lambda)}(y^*) - F(y^*; v_n^{(\lambda)}(t^*,\cdot))| 
 \le C(1+\lambda + \lambda^2)R_nh_n. 
\] 
Therefore, by \eqref{eq:3.3}, 
\begin{align*}
 &|-\partial_tv_n^{(\lambda)}(t_n,x_n) + F(t_n,x_n; v_n^{(\lambda)}(t_n,\cdot))| \\ 
 &\le |-\partial_tv_n^{(\lambda)}(y_n) + F(y_n; v_n^{(\lambda)}(t_n,\cdot)) + \partial_tv_n^{(\lambda)}(y) - F(y; v_n^{(\lambda)}(t,\cdot))| \\ 
 &\quad + |-\partial_tv_n^{(\lambda)}(y) + F(y; v_n^{(\lambda)}(t,\cdot)) + \partial_tv_n^{(\lambda)}(y^*) - F(y^*; v_n^{(\lambda)}(t^*,\cdot))| \\
 &\quad + |-\partial_tv_n^{(\lambda)}(y^*) + F(y^*; v_n^{(\lambda)}(t^*,\cdot))| \\ 
 &\le C(1+\lambda + \lambda^2)(|y-y_n| + R_nh_n) + \gamma_n^{(\lambda)}, 
\end{align*}
where 
\[
 \gamma_n^{(\lambda)}=\max_{\ell=1,\ldots,M}\left|-\partial_tv_n^{(\lambda)}(\tilde{y}^{\ell})+ F(\tilde{y}^{\ell}; v_n^{(\lambda)}(\tilde{t}^{\ell},\cdot))\right|
 \vee|v_n^{(\lambda)}(T,\tilde{x}^{\ell}) - g(\tilde{x}^{\ell})|.  
\]

We will show that $\gamma_n=\gamma_n^{(\lambda)}\to 0$, with suitable $\lambda$.   
Denote by $I[v](y)$ the interpolation of $v$ on $\Gamma_n$ with the reproducing 
kernel $\Phi$, i.e., 
\[
 I[f](y) = \sum_{j=1}^N(K^{-1}\tilde{f})^j\Phi(y-y^j), \quad y\in\mathbb{R}^{d+1}, 
\]
where $\tilde{f}=(f(y^1),\ldots,f(y^N))^{\mathsf{T}}$ for $f:\mathbb{R}^{d+1}\to\mathbb{R}$.  
Then, by Lemma 3.2 below, 
\begin{equation}
\label{eq:3.4}
 \sup_{y\in Q_n}|D^{\alpha}v(y) - D^{\alpha}I[v](y)|\le C(R_nh_n)^{\tau-(d+1)/2 - |\alpha|_1}\|v\|_{\mathcal{H}}, \quad |\alpha|_1\le 2.  
\end{equation}
It should be emphasized that the constant $C$ in \eqref{eq:3.4} does not depend on $n$. 

Put $\theta=K^{-1}\tilde{v}$. Then, by \eqref{eq:3.4} and the interpolation property, $u(t,x):=I[v](t,x)$ satisfies
\begin{align*}
 \left| - \partial_t u(\tilde{t}^{\ell},\tilde{x}^{\ell}) + F(\tilde{t}^{\ell},\tilde{x}^{\ell}; u(\tilde{t}^{\ell}))\right| 
  &= \left| - \partial_t u(\tilde{t}^{\ell},\tilde{x}^{\ell}) + F(\tilde{t}^{\ell},\tilde{x}^{\ell}; u(\tilde{t}^{\ell})) 
       + \partial_t v(\tilde{t}^{\ell},\tilde{x}^{\ell}) - F(\tilde{t}^{\ell},\tilde{x}^{\ell}; v(\tilde{t}^{\ell}))\right| \\ 
 &\le C\sum_{|\alpha|_1\le 2}\left| D^{\alpha}v(\tilde{t}^{\ell},\tilde{x}^{\ell}) - D^{\alpha}I[v](\tilde{t}^{\ell},\tilde{x}^{\ell})\right|  \\ 
 &\le C^{\prime} (R_nh_n)^{\tau-(d+1)/2-2}\|v\|_{\mathcal{H}}
\end{align*}
and 
\[
 |u(T,\tilde{x}^{\ell}) - g(\tilde{x}^{\ell})|\le C^{\prime} (R_nh_n)^{\tau-(d+1)/2}\|v\|_{\mathcal{H}}
\]
for some constant $C^{\prime}>0$. By the assumption imposed in the theorem and Corollary 10.25 in \cite{wen:2010}, there exists a positive constant $C^{\prime\prime}$ 
such that  
\[
 \theta^{\mathsf{T}}K\theta = \|u\|_{\mathcal{H}}^2\le \|v\|_{\mathcal{H}}^2\le C^{\prime\prime}
\]
for $s,t\in [0,T]$. 
Thus $(\theta, \gamma)\in \mathbb{R}^N\times [0,\infty)$ satisfies the constraint in (P$_\lambda$) with the choice $\lambda=\sqrt{C^{\prime\prime}}$, 
where 
\[
 \gamma=\max_{\ell=1,\ldots,M}\left| - \partial_t u(\tilde{t}^{\ell},\tilde{x}^{\ell}) + F(\tilde{t}^{\ell},\tilde{x}^{\ell}; u(\tilde{t}^{\ell}))\right|
  \vee |u(T,\tilde{x}^{\ell}) - g(\tilde{x}^{\ell})|. 
\]
This means that $\gamma_n^{(\lambda)}\le C^{\prime}h_n^{\tau-(d+1)/2-2}\|v\|_{\mathcal{H}} + \varepsilon_n\to 0$, as $n\to\infty$. 
We conclude that the condition (B1) is satisfied for $(t,x)\in (0,T)\times\mathbb{R}^d$. 

Step (ii). 
We will check (B1) for $t=0$. Let $x\in\mathbb{R}^d$ and $n_0\in\mathbb{N}$ such that $|x-\bar{x}|<R_n$ for $n\ge n_0$. 
Let $(t_n,x_n)\in [0,T)\times\mathbb{R}^d$ be such that $(t_n,x_n)\to y:=(0,x)$. 
Take $(t^*,x^*)\in Q_n$ such that $|y-y^*|=\min_{j=1,\ldots,N_n}|y-y^j|$. 
Then, since $\sup_{y^{\prime}\in Q_n}\min_{j=1,\ldots,N_n}|y^{\prime}-y^j| = \sup_{y^{\prime}\in\overline{Q_n}}\min_{j=1,\ldots,N_n}|y^{\prime}-y^j|$,  
as in Step (i) we get 
\begin{align*}
 &|-\partial_tv_n^{(\lambda)}(t_n,x_n) + F(t_n,x_n; v_n^{(\lambda)}(t_n,\cdot))| \\ 
 &\le |-\partial_tv_n^{(\lambda)}(t_n,x_n) + F(t_n,x_n; v_n^{(\lambda)}(t_n,\cdot)) + \partial_tv_n^{(\lambda)}(0,x) - F(0,x; v_n^{(\lambda)}(0,\cdot))| \\ 
 &\quad + |-\partial_tv_n^{(\lambda)}(0,x) + F(0,x; v_n^{(\lambda)}(0,\cdot)) + \partial_tv_n^{(\lambda)}(t^*,x^*) - F(t^*,x^*; v_n^{(\lambda)}(0,\cdot))| \\
 &\quad + |-\partial_tv_n^{(\lambda)}(t^*,x^*) + F(t^*,x^*; v_n^{(\lambda)}(0,\cdot))| \\ 
 &\le C(1+\lambda + \lambda^2)(t_n+ |x-x_n| +R_n h_n) + \gamma_n^{(\lambda)}, 
\end{align*}
leading to the claim.

Step (iii). 
Finally we will show that the consistency (B2) holds. Let $x\in\mathbb{R}^d$ be fixed, and $(t_n,x_n)$ be such that $(t_n, x_n)\to (T,x)$, as $n\to\infty$. 
There exists $n_0\ge 1$ such that $|x-\bar{x}|<R_n$ for any $n\ge n_0$. Then, again by \eqref{eq:3.4}, 
\begin{align*}
 \left|\sum_{j=1}^N\theta_j\Phi((t_n,x_n)-y^j) - g(x_n)\right| 
 &\le |I[v(t_n)](x_n) - v(t_n,x_n)| + |v(t_n,x_n) - g(x_n)| \\
 &\le C\lambda (R_nh_n)^{\tau-(d+1)/2} + C\lambda (T-t_n)
 \to 0,
\end{align*}
as $n\to\infty$. Thus (B2) follows. 
\end{proof}

We have used the following result in the proof of Theorem \ref{thm:3.1}: 
\begin{lem}
\label{lem:3.2}
There exists $n_0\in\mathbb{N}$ such that 
for any multi-index $\alpha$ with $|\alpha|_1\le 2$ and $f\in\mathcal{H}$, 
we have 
\begin{equation*}
 |D^{\alpha}f(t,x)-D^{\alpha}I[f](t,x)|\le 
 C(R_nh_n)^{\tau- (d+1)/2 -|\alpha|_1}\|f\|_{\mathcal{H}}, 
 \quad (t,x)\in \overline{Q_n}, \;\; n\ge n_0.  
\end{equation*}
\end{lem}
\begin{proof}
This result is reported in \cite[Corollary 11.33]{wen:2010} for more general domains. 
However, a simple application of that result leads to an ambiguity of the dependence of 
the constant $C$ on $R$. 
Here we will confirm that we can take $C$ to be independent of $R$. 
We prove the claim by a minor modification of the proof of Lemma 3.6 in Nakano \cite{nak:2019a}.  

Let $f\in H^{\tau}((0,T)\times\mathbb{R}^d)$ with $f|_{\Gamma_n}=0$. 
Define $\tilde{f}(y)=f(\iota(y))$, $y\in Q$. 
Then, $\tilde{f}|_{\Gamma}=f|_{\Gamma_n}=0$.  
Since $h_n\to 0$ as $n\to\infty$ and $\tau> 3+(d+1)/2$, 
we can apply \cite[Theorem 11.32]{wen:2010} to 
$\tilde{f}$ to obtain 
\begin{equation}
\label{eq:3.2lem}
 |D^{\alpha}\tilde{f}(y)|\le Ch_n^{\tau-(d+1)/2-|\alpha|_1}
 |\tilde{f}|_{H^{\tau}(Q)}, \quad y\in Q, \;\; n\ge n_0
\end{equation}
for some $n_0\in\mathbb{N}$.
It is straightforward to see that for $\alpha=(\alpha_0,\alpha_1,\ldots,\alpha_d)$,
\[
 D^{\alpha}f(y)=\frac{\partial^{|\alpha|_1}f}{\partial_t^{\alpha_0}\partial_{x_1}^{\alpha_1}\cdots\partial_{x_d}^{\alpha_d}}(y)
 = \left(\frac{T}{2}\right)^{-\alpha_0}R^{-(\alpha_1+\cdots +\alpha_d)}(D^{\alpha}\tilde{f})(\iota(y)).
\]
Substituting these relations into \eqref{eq:3.2lem}, we have
\begin{equation}
\label{eq:3.3lem}
 |D^{\alpha}f(y)|\le CR^{-|\alpha|_1}h_n^{\tau-(d+1)/2- |\alpha|_1}|\tilde{f}|_{H^{\tau}(Q)}
 \le C(R_nh_n)^{\tau-(d+1)/2-|\alpha|_1}\|f\|_{H^{\tau}(Q_n)}, 
 \quad y\in Q_n. 
\end{equation}
This and \cite[Corollary 10.25]{wen:2010} yield 
\begin{align*}
 |D^{\alpha}f(y)-D^{\alpha}I[f](y)|
 &\le C(R_nh_n)^{\tau-(d+1)/2-|\alpha|_1}\|f-I(f)\|_{\mathcal{H}} \\ 
 &\le C(R_nh_n)^{\tau-(d+1)/2-|\alpha|_1}\|f\|_{\mathcal{H}}
\end{align*}
for any $y\in\overline{Q_n}$. Thus the lemma follows. 
\end{proof}


\subsection{Convergence for Hamilton--Jacobi--Bellman equations}\label{sec:3.2}

In this section, we apply the kernel-based approximation method introduced in the previous section to Hamilton--Jacobi--Bellman (HJB) equations. 
By applying the result of Section \ref{sec:2.2}, we obtain not only convergence but also quantitative error estimates. 

Here we consider HJB equation
\begin{equation}
\label{eq:3.5}
\left\{
\begin{split}
 & \partial_t v + \inf_{a\in A}\mathcal{L}^av(t,x)=0, \quad (t,x)\in
 [0,T)\times \mathbb{R}^d, \\
 &v(T,x)=g(x), \quad x\in\mathbb{R}^d,
\end{split}
\right.
\end{equation}
where $A\subset \mathbb{R}^{d_1}$ is compact, and for $a\in A$,  
\[
 \mathcal{L}^a\phi(t,x)=b(t,x,a)^{\mathsf{T}}D_x\phi(t,x) + \frac{1}{2}\mathrm{tr}(\sigma\sigma^{\mathsf{T}}(t,x,a)D^2_{xx}\phi(t,x)) + f(t,x,a), \quad (t,x)\in [0,T)\times\mathbb{R}^d,
\] 
with some continuous functions $b:[0,T]\times\mathbb{R}^d\times A\to\mathbb{R}^d$,
$\sigma: [0,T]\times\mathbb{R}^d\times A\to\mathbb{R}^{d\times d_2}$, and $f:[0,T]\times\mathbb{R}^d\times A\to\mathbb{R}$.

Let $\Phi$, $Q_n$ and $\Gamma_n$ be as in Section \ref{sec:3.1}. Here, the parameter $\lambda$ in Section \ref{sec:3.1} is also assumed to depend on $n$, 
and so we write $\lambda=\lambda_n$. Thus the problem ($P_{\lambda}$) becomes 
\begin{equation}
\tag{$P_n$}\label{eq:Pn}
\left\{
\begin{aligned}
 &\text{minimize}\;\; \gamma, \\
 &\text{s.t.}\;\;  (\theta, \gamma) \in \mathbb{R}^N\times [0,\infty), \;\; \theta^{\mathsf{T}}K\theta\le \lambda^2, \;\;
  -\gamma \mathbf{1}_N\le \tilde{K}\theta - \tilde{g} \le \gamma \mathbf{1}_N, \;\;
  -\gamma \mathbf{1}_N\le \mathbf{F}(\theta) \le \gamma \mathbf{1}_N,
\end{aligned}
\right.
\end{equation}
where 
\[
 F(t,x; \phi(t))= - \inf_{a\in A}\mathcal{L}^a\phi(t,x). 
\]

Our candidate of approximate solution is now constructed as follows: 
given $\lambda=\lambda_n>0$ and $\varepsilon=\varepsilon_n>0$, 
take $(\theta^{(n)}, \gamma^{(n)})\in\mathbb{R}^N\times \mathbb{R}$ to be 
an $\varepsilon_n$-optimal solution of ($P_n$). Then, define $v_n\in\mathcal{H}$ by 
\[
 v_n(y) = \sum_{j=1}^N\theta_j^{(n)}\Phi(y-\xi^j), \quad y\in [0,T]\times\mathbb{R}^d.
\] 

We shall give an error estimation of the methods above using the result in Section \ref{sec:2.2}. 
\begin{thm}
\label{thm:3.3}
Suppose that \rm{(A4)} and \rm{(A6)} hold and that there exists a classical solution $v$ of $\eqref{eq:3.5}$ and $\beta_0\in (0,1]$ such that
for $z=\partial_tv, D_xv, D_x^2v$, 
\[
\sup_{(t,x), (s,y)\in [0,T]\times \mathbb{R}^d}\frac{|z(t,x) - z(s,y)|}{|t-s|^{\beta_0/2} + |x-y|^{\beta_0}} < \infty.
\]
Then, for any compact set $Q\subset [0,T]\times\mathbb{R}^d$ there exist a constant $C>0$ and $n_0\in\mathbb{N}$ such that for $n\ge n_0$, 
\begin{equation}
\label{eq:3.8}
\begin{aligned}
 &\sup_{(t,x)\in Q}|v(t,x)-v_n(t,x)| \\ 
 & \le (1+T)\varepsilon_n + C\lambda_n (R_nh_n)^{\tau - (d+1)/2-2} + C(\lambda_nR_n^{\tau-d/2})^{-\beta_0/(2\tau-1)}
    + C\lambda_nR_nh_n + C \lambda_ne^{- R_n^2/(8C_1^2T)}.
\end{aligned}
\end{equation}
\end{thm}
\begin{rem}
From \eqref{eq:3.8} a set of sufficient conditions for the convergence is that $\lambda_n\to\infty$, 
$R_n\to\infty$, $\lambda_nR_nh_n\to 0$ and $\lambda_ne^{- R_n^2/(8C_1^2T)}\to 0$ as $n\to\infty$.
\end{rem}
\begin{proof}[Proof of Theorem \rm{\ref{thm:3.3}}]
Let $\Psi\in C^{\infty}_0(\mathbb{R}\times\mathbb{R}^d)$ supported in $[0,1]\times\{x\in\mathbb{R}^d: |x|\le 1\}$, with unit integral; this is a mollifier and is unrelated to the radial basis function $\Phi$ of Section~\ref{sec:3.1}.
Define $\Psi_n$ by $\Psi_n(t,x)=\rho^{d+2}\Psi(\rho^2 t, \rho x)$, $(t,x)\in \mathbb{R}\times\mathbb{R}^d$, where $\rho=\rho_n$ will be specified later.
Then, consider the convolution $w_n$ of $\Psi_n$ with an extension of $v$ to $[-1/\rho^2, T]$, defined by
\[
 w_n(t,x):=\int_{\mathbb{R}}\int_{\mathbb{R}^d}v(t-s,x-y)\Psi_n(s,y)\,ds\,dy, \quad (t,x)\in\mathbb{R}\times\mathbb{R}^d.
\]
Note we can always take an extension of $v(t,x)$ to $[-1/\rho^2,T]\times\mathbb{R}^d$, which is still in $C^{1+\beta_0/2, 2+\beta_0}([-1/\rho^2, T]\times\mathbb{R}^d)$.
It is straightforward to see that 
\[
 \|w_n\|_{C^{\tau}([0,T]\times\mathbb{R}^d)}\le C\rho^{2\tau-1},
\]
and 
\[
 \|v - w_n\|_{C^{1,2}([0,T]\times\mathbb{R}^d)}\le C\rho^{-\beta_0}. 
\]
Let $\zeta\in C^{\infty}_0(\mathbb{R}^{d+1})$ such that $0\le \zeta\le 1$ on $\mathbb{R}^{d+1}$, $\zeta=1$ on $Q$, and 
$\zeta=0$ on $\{y\in\mathbb{R}^{d+1}: |y|_{\infty}>2\}$. 
Consider $\tilde{w}_n(t,x) := w_n(t,x)\zeta(\iota^{-1}(t,x))$, $(t,x)\in\mathbb{R}^{d+1}$. 
Observe $|D^{\alpha}\tilde{w}_n(t,x)|\le C\rho^{2\tau-1}R^{-\tau}$ for $|\alpha|_1\le\tau$, whence 
\[
 \|\tilde{w}_n\|_{\mathcal{H}}\le C\rho^{2\tau-1}R^{-\tau+d/2}. 
\]
In view of this, we determine $\rho>0$ by $\rho^{2\tau-1}=C\lambda_n R_n^{\tau-d/2}$. 
Put $u_n=I[w_n]=I[\tilde{w}_n]$. 
Then, by Corollary 10.25 in \cite{wen:2010}, 
\begin{equation}
\label{eq:3.9}
 \|u_n\|_{\mathcal{H}}\le \lambda_n. 
\end{equation}
Further, by Lemma \ref{lem:3.2}, for $|\alpha|_1\le 2$ and $y\in\overline{Q}_n$, 
\begin{align*}
|D^{\alpha}w_n(y) - D^{\alpha}u_n(y)|&= |D^{\alpha}\tilde{w}_n(y) - D^{\alpha}I[\tilde{w}_n](y)| \\ 
&\le C(R_nh_n)^{\tau-(d+1)/2 - 2}\|\tilde{w}_n\|_{\mathcal{H}} \\ 
&\le C\lambda_n (R_nh_n)^{\tau-(d+1)/2 - 2}. 
\end{align*}
Using this, we get 
\begin{equation}
\label{eq:3.10}
\begin{aligned}
 |u_n(T,x) - g(x)|&\le |u_n(T,x) - w_n(T,x)| + |w_n(T,x)-g(x)| \\ 
 &\le C\lambda_n (R_nh_n)^{\tau - (d+1)/2} + C\rho^{-\beta_0} \\ 
 &\le C\lambda_n (R_nh_n)^{\tau - (d+1)/2} + C(\lambda_nR_n^{\tau-d/2})^{-\beta_0/(2\tau-1)}
\end{aligned}
\end{equation}
Further, 
\begin{equation}
\label{eq:3.11}
\begin{aligned}
 |(\partial_t + \inf_{a\in A}\mathcal{L}^a)u_n(t,x)| &\le C\|u_n- w_n\|_{C^{1,2}([0,T]\times\mathbb{R}^d)} + C\|v - w_n\|_{C^{1,2}([0,T]\times\mathbb{R}^d)} \\ 
 &\le C\lambda_n (R_nh_n)^{\tau - (d+1)/2-2} + C\rho^{-\beta_0} \\ 
 &\le C\lambda_n (R_nh_n)^{\tau - (d+1)/2-2} + C(\lambda_nR_n^{\tau-d/2})^{-\beta_0/(2\tau-1)}. 
\end{aligned}
\end{equation}
Therefore, \eqref{eq:3.9}--\eqref{eq:3.11} means that if $\tilde{\theta}\in\mathbb{R}^N$ is such that $u_n(y)=\sum_{j=1}^N\tilde{\theta}_j\Phi(y-y^j)$ 
and 
\[ 
 \tilde{\gamma}=\max_{\ell=1,\ldots,M}|(\partial_t + \inf_{a\in A}\mathcal{L}^a)u_n(\tilde{t}^{\ell},\tilde{x}^{\ell})|\vee |u_n(T,\tilde{x}^{\ell}) - g(\tilde{x}^{\ell})|
\]
then $(\tilde{\theta},\tilde{\gamma})$ satisfies the constraints in ($P_n$), and so 
\[
\gamma^{(n)}\le \varepsilon_n + C\lambda_n (R_nh_n)^{\tau - (d+1)/2-2} + C(\lambda_nR_n^{\tau-d/2})^{-\beta_0/(2\tau-1)}.
\]
For any $y=(t,x)\in Q_n$ take a nearest $\tilde{y}\in\Gamma_n$ such that $\min_{j=1,\ldots,N}|y-y^j|=|y-\tilde{y}|$. Then 
\begin{align*}
 |(\partial_t + \inf_{a\in A}\mathcal{L}^a)v_n(y)|\vee |v_n(T,x) - g(x)| 
 &\le \gamma^{(n)} + C\|v_n\|_{\mathcal{H}}|y-\tilde{y}|\le \gamma^{(n)} + C\lambda_nR_nh_n. 
\end{align*}
This leads to the claim. 
\end{proof}

\section{Numerical experiments}\label{sec:4}

This section provides a minimal numerical validation of the proposed method.
Our goal is not to claim state-of-the-art speed or accuracy, but to confirm that
the non-monotone approximation is computationally feasible and consistent with the theoretical convergence framework.

Because the method requires solving a large-scale constrained optimization problem,
we report not only accuracy but also computational cost and optimization health
(constraint violation and final objective value).

In what follows, the function $\Phi$ is 
assumed to be a Wendland kernel $\Phi(y)=\phi(|y|)$ (see \cite[Chapter 9]{wen:2010}).

\subsection{Matrix representation}\label{subsec:matrix_rep}

First let us give matrix representations of the derivative of the kernel that is easier to handle numerically.  
Consider the function $\phi^{(1)}(r):=\phi^{\prime}(r)/r$, $r\ge 0$. 
By definition of $\phi$, the function $\phi^{(1)}$ is continuous on 
$[0,\infty)$ and supported in $[0,1]$. With this function, we have 
\[
 \partial_{y_k}\Phi(y)=\phi^{(1)}(|y|)y_k, \quad y=(y_0,y_1,\ldots,y_d)\in\mathbb{R}^{d+1}, 
\]
where $y_0$ corresponds to the time variable $t$ and $y_i$ corresponds to the $i$-th element of the spatial variable $x$ for $i=1,\ldots,d$. 
Thus, 
\[
 B_k:= \left((\partial_{y_k}\Phi)(\tilde{y}^{\ell} - y^j)\right)_{{1\le \ell\le M}\atop{1\le j\le N}}
  = \tilde{G}_{k}K_1 - K_1G_{k}, 
\]
where $K_1 = \{\phi^{(1)}(|\tilde{y}^{\ell} - y^j|)\}_{{1\le\ell\le M}\atop{1\le j\le N}}$, 
$G_k = \mathrm{diag}(y^{1}_k,\ldots,y^{N}_k)$
and $\tilde{G}_{k} = \mathrm{diag}(\tilde{y}^{1}_{k},\ldots,\tilde{y}^{N}_{k})$ for $k=0,1,\ldots,d$. 
Similarly, 
\[
 \partial^2_{y_m y_{k}}\Phi(y) = 
 \begin{cases}
  \phi^{(1)}(|y|) + \phi^{(2)}(|y|)y_m^2, & (k=m), \\ 
  \phi^{(2)}(|y|)y_m y_{k}, & (k \neq m), 
 \end{cases} 
\]
where 
\[
 \phi^{(2)}(r) = \frac{1}{r}\frac{d\phi^{(1)}}{dr}(r), \quad r\ge 0.
\]
Notice that $\phi^{(2)}$ is also continuous on $[0,\infty)$ and supported in $[0,1]$. 
Thus, 
\[
 B_{m k}:=\left((\partial^2_{y_my _{k}}\Phi)(\tilde{y}^{\ell}-y^{j})\right)_{{1\le\ell\le M}\atop{1\le j\le N}}
\]
is given by 
\[
 B_{mm} = K_1 + \tilde{G}_m^2K_2 - 2\tilde{G}_mK_2G_m + K_2G_m^2
\]
and for $m\neq k$, 
\[
 B_{m k}
 = \tilde{G}_m\tilde{G}_{k}K_2 - \tilde{G}_mK_2G_{k}-\tilde{G}_{k}K_2G_m + K_2G_mG_{k}
\]
with $K_2 =\{\phi^{(2)}(|\tilde{y}^{\ell}-y^{j}|)\}_{{1\le\ell\le M}\atop{1\le j\le N}}$.

\subsection{Test problem and setup}
\label{subsec:test_setup}

We solve the problem ($P_n$) for the following $d$-dimensional equation,
adopted from \cite{guo-etal:2015}: for $d\in\mathbb{N}$,
\begin{equation}
\label{eq:4.testHJB}
\left\{
\begin{aligned}
 &-\partial_t v - \frac{1}{2}\sup_{0\le\sigma\le 1/5}
   \mathrm{tr}\bigl(\sigma^2\,D^2_{xx}v\bigr)
   + G(v,D_xv)=0,\quad (t,x)\in[0,1)\times\mathbb{R}^d,\\
 &v(1,x)=\sin\Bigl(1+\textstyle\sum_{i=1}^{d}x_i\Bigr),\quad x\in\mathbb{R}^d,
\end{aligned}
\right.
\end{equation}
where
\(
 G(z,p)=\tfrac{1}{d}\sum_{i=1}^{d}p_i-\tfrac{d}{2}\inf_{0\le\sigma\le 1/5}(\sigma^2 z)
\)
for $z\in\mathbb{R}$, $p=(p_1,\ldots,p_d)\in\mathbb{R}^d$.
The unique classical solution is
\(
 v(t,x)=\sin\bigl(t+\textstyle\sum_{i=1}^{d}x_i\bigr).
\)

We work out the case $d=1$ in detail in the present and the next subsection;
the $d=2$ extension is reported in Section~\ref{subsec:d2}.

For $d=1$, the function $\Phi$ is taken as the Wendland radial basis function
$\Phi(y)=\phi(|y|)$, $y\in\mathbb{R}^2$, with
\[
 \phi(r)=(1-r/\varepsilon)_+^5\bigl(8(r/\varepsilon)^2+5r/\varepsilon+1\bigr),\quad r\ge 0,
\]
for a given bandwidth $\varepsilon>0$.
$\Phi$ is positive definite on $\mathbb{R}^2$ and belongs to $C^4(\mathbb{R}^2)$,
as required by the abstract framework of Section~\ref{sec:3}.

The collocation set $\Gamma_n$ is given by a uniform tensor-product grid where
$N^t+1$ equally spaced time nodes $t_k = k/N^t$
($k=0,\ldots,N^t$) are combined with $N^x$ equally spaced points
per spatial axis on $[-R,R]$, giving $N^x$ spatial nodes and
a total of $N = (N^t+1)\cdot N^x$ collocation points.
An additional set $\Gamma^{\rm extra}$ of $N^{\rm extra}$ points is drawn from a
$2$-dimensional Sobol' sequence (with skip $10^3$ and leap $10^2$)
and mapped to $[0, 1-\delta_0)\times[-R,R]$, where
$\delta_0 = 10^{-8}$ ensures that the terminal time $t=1$ is
excluded.
These points are used to impose supplementary PDE residual
constraints that are not tied to the collocation grid.
Then we use $\Gamma^e=\Gamma\cup\Gamma^{\rm extra}$ for PDE evaluation, error assessment
and the early stopping criterion. 

For practical solution of ($P_n$) we adopt a smooth penalty reformulation:
the constraint $\theta^{\mathsf{T}}K\theta\le \lambda^2$ is replaced by a quadratic penalty,
the inequality bounds on the residual and terminal mismatch are absorbed into a
least-squares loss, and the resulting unconstrained problem
\begin{equation}
\label{eq:Pn-penalty}
 \min_{\theta\in\mathbb{R}^N}\;
 \frac{1}{M}\sum_{\ell=1}^M|\mathbf{F}^{\ell}(\theta)|^2
   + \frac{1}{M}\sum_{\ell=1}^M|\tilde K_{\ell:}\theta-\tilde g_{\ell}|^2
   + \mu\,\bigl(\theta^{\mathsf{T}}K\theta-\lambda^2\bigr)_+^2
\end{equation}
is solved by the limited-memory BFGS method \cite{liu-noc:1989},
with an initial point produced by the collocation method of \cite{kan:1990b,nak:2017}.

The penalty weight $\mu$ is chosen as $\mu=10^2$.
A sensitivity analysis on the $N=1056$ instance (S5 in Table~\ref{tab:main_result})
shows that $E_{\mathrm{RMS}}$ varies by less than $5\%$ across $\mu\in\{10^0,10^1,\ldots,10^5\}$,
with the constraint $\theta^{\mathsf{T}}K\theta\le \lambda^2$ satisfied to within $10^{-4}$
for all tested $\mu\ge 10$;
the choice $\mu=10^2$ is the smallest value for which the constraint is satisfied to
machine precision while the optimizer converges in a few hundred L-BFGS iterations.

Once the optimizer returns $\hat\theta$, the residual bound $\gamma_n$ is computed
\emph{a posteriori} as
\(
\gamma_n=\max_{\ell}\bigl(|\mathbf{F}^{\ell}(\hat\theta)|\vee|\tilde K_{\ell:}\hat\theta-\tilde g_{\ell}|\bigr).
\)

The main parameters are the number $N$ of space-time collocation points,
the radius $R$ in the space domain, and the bandwidth $\varepsilon$ of the kernel.
In all experiments below we set $R=3$, $\lambda=2R=6$, and $\varepsilon=0.5$.
Eight settings (S1--S8), distinguished by the triple $(N^t, N^x, N^{\rm extra})$, are reported in Table~\ref{tab:main_result} of Section~\ref{subsec:results}.

\subsection{Numerical results in $d=1$}
\label{subsec:results}

We report the errors
\begin{align*}
E_{\infty} &:= \max_{\ell=1,\ldots,M}|v_n(\tilde{t}^{\ell},\tilde{x}^{\ell})-v(\tilde{t}^{\ell},\tilde{x}^{\ell})|,\\
E_{\mathrm{RMS}} &:= \left(\frac{1}{M}\sum_{\ell=1}^M|v_n(\tilde{t}^{\ell},\tilde{x}^{\ell})-v(\tilde{t}^{\ell},\tilde{x}^{\ell})|^2\right)^{1/2},
\end{align*}
together with the residual sup-norm $\mathrm{Res}_{\infty}:=\max_{\ell}|\mathbf{F}^{\ell}(\hat\theta)|$,
the residual bound $\gamma_n$ defined after \eqref{eq:Pn-penalty},
the iteration count, and the total CPU time.

Table~\ref{tab:main_result} presents the main results.
The table is automatically generated from the output of the implementation
that accompanies this paper.
All CPU times reported in Sections~\ref{subsec:results}--\ref{subsec:d2}
were measured on a single core of an Apple~M3~Pro chip
with 18\,GB of unified memory.

\begin{table*}[tb]
  \centering
  \caption{Main results for the test problem of Section~\ref{subsec:test_setup}: approximation errors, residual bound $\gamma_n$, iteration count, and CPU time.}
  \label{tab:main_result}
  \begin{tabular}{ccccccccccc}
    \toprule
    Case & $N^t$ & $N^x$ & $N^{\mathrm{extra}}$ & $N$ & $\gamma_n$ & $E_{\mathrm{RMS}}$ & $E_{\infty}$ & $\mathrm{Res}_{\infty}$ & Iter. & CPU [s] \\
    \midrule
    S1 & 32 & 16 & 100 & 528 & 2.043e+00 & 4.365e-01 & 8.561e-01 & 2.043e+00 & 5000 & 7.6 \\
    S2 & 16 & 16 & 500 & 272 & 1.578e+00 & 5.970e-01 & 9.438e-01 & 1.578e+00 & 3037 & 2.9 \\
    S3 & 16 & 32 & 100 & 544 & 5.771e-01 & 1.438e-01 & 6.343e-01 & 5.771e-01 & 598 & 1.3 \\
    S4 & 16 & 32 & 500 & 544 & 7.256e-01 & 1.469e-01 & 6.477e-01 & 7.256e-01 & 641 & 2.6 \\
    S5 & 32 & 32 & 100 & 1056 & 8.253e-01 & 1.419e-01 & 6.072e-01 & 8.253e-01 & 1047 & 10.0 \\
    S6 & 16 & 64 & 100 & 1088 & 5.333e-01 & 9.707e-02 & 4.119e-01 & 5.333e-01 & 419 & 4.1 \\
    S7 & 16 & 64 & 500 & 1088 & 6.203e-01 & 1.122e-01 & 4.953e-01 & 6.203e-01 & 674 & 8.8 \\
    S8 & 32 & 64 & 100 & 2112 & 8.608e-01 & 8.804e-02 & 3.431e-01 & 8.608e-01 & 308 & 11.8 \\
    \bottomrule
  \end{tabular}
\end{table*}

\begin{figure}[h]
\centering
\begin{minipage}[b]{0.49\columnwidth}
  \centering
  \includegraphics[width=0.95\linewidth]{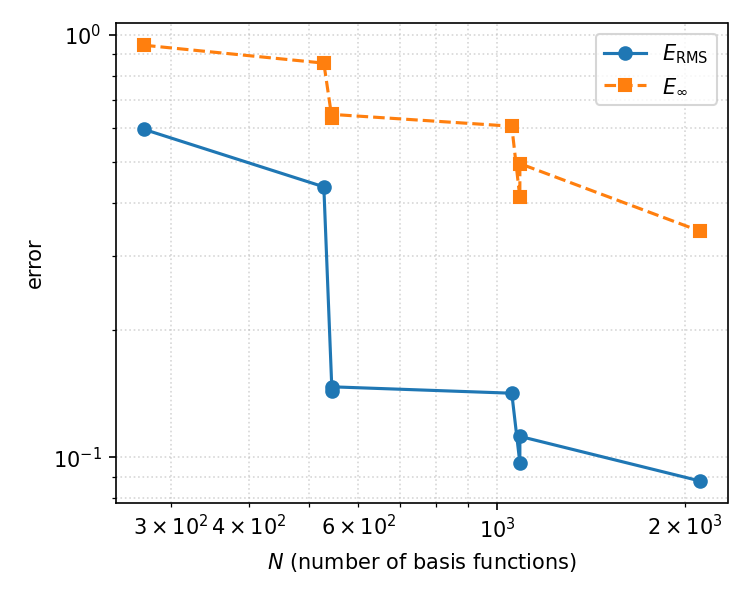}
  \caption{The errors $E_{\mathrm{RMS}}$ and $E_{\infty}$ versus the number of basis functions $N$.}
  \label{fig:error_vs_size}
\end{minipage}\hfill
\begin{minipage}[b]{0.49\columnwidth}
  \centering
  \includegraphics[width=0.95\linewidth]{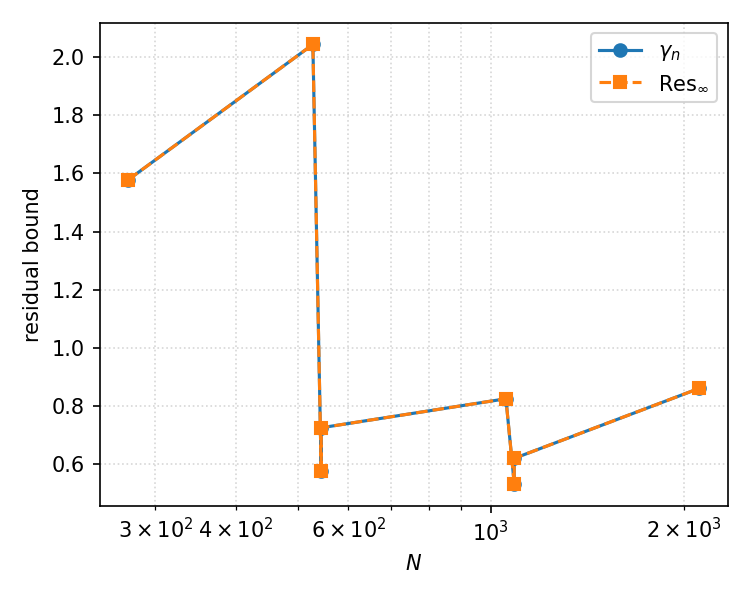}
  \caption{The residual bound $\gamma_n$ and $\mathrm{Res}_{\infty}$ versus $N$.}
  \label{fig: gamma_vs_size}
\end{minipage}
\end{figure}

\begin{figure}[h]
\centering
\begin{minipage}[b]{0.49\columnwidth}
  \centering
  \includegraphics[width=0.95\linewidth]{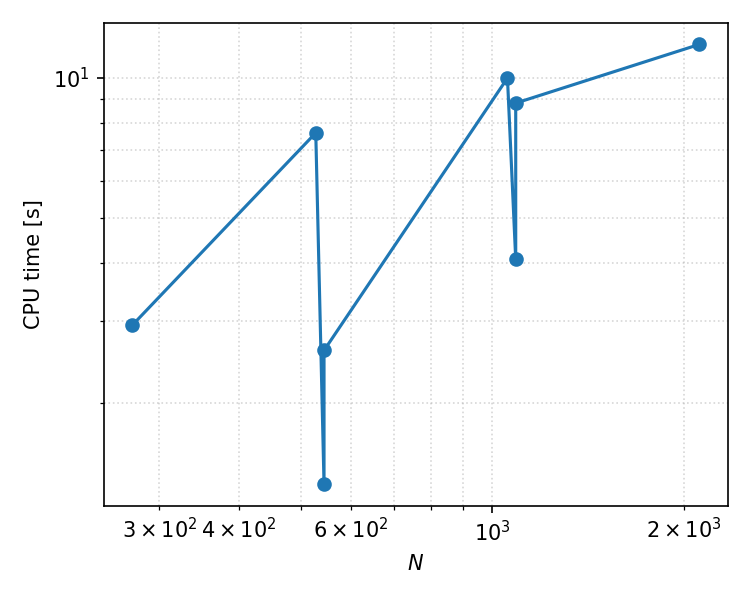}
  \caption{CPU time versus $N$.}
  \label{fig:time_vs_size}
\end{minipage}
\end{figure}

The data of Table~\ref{tab:main_result} support the qualitative behaviour predicted
by Theorem~\ref{thm:3.1} and the rate analysis of Theorem~\ref{thm:3.3}.

\smallskip
\noindent\emph{Decrease of errors with $N$.}
Both $E_{\mathrm{RMS}}$ and $E_{\infty}$ decrease as the number of basis functions $N$ grows
(Figure~\ref{fig:error_vs_size}): from $E_{\mathrm{RMS}}\approx 6\times 10^{-1}$
at $N=272$ down to $E_{\mathrm{RMS}}\approx 9\times 10^{-2}$ at $N=2112$, with a similar
trend for $E_{\infty}$.
The decrease is consistent with the predicted dependence on $N$ in Theorem~\ref{thm:3.1};
saturating below $10^{-1}$ at $N\approx 10^{3}$ reflects the bias introduced by the
fixed kernel width $\varepsilon=0.5$ relative to the wavelength of the oscillatory target $\sin(t+x)$.

\smallskip
\noindent\emph{Decrease of the residual bound.}
The residual bound $\gamma_n$ also decreases overall with $N$
(Figure~\ref{fig: gamma_vs_size}), from $\gamma_n\approx 1.6$ at $N=272$ down to
$\gamma_n\approx 0.55$ at $N\ge 1000$. The non-monotone behaviour at the smallest scales
reflects the role of the supplementary collocation set $\Gamma^{\rm extra}$:
a larger $\Gamma^{\rm extra}$ provides a sharper estimate of the residual sup-norm but
does not change the basis dimension, hence the slight increase of $\gamma_n$ at fixed $N$
when $|\Gamma^{\rm extra}|$ grows (cf.\ S3 vs S4 and S6 vs S7).

\smallskip
\noindent\emph{Computational cost.}
The unconstrained penalty formulation \eqref{eq:Pn-penalty} solves entirely on a single CPU
core in a few seconds per case (Figure~\ref{fig:time_vs_size}); the dominant cost is the
construction of the kernel matrix and its derivative tables, which scales as $O(MN)$.
The improvement over a constrained sequential quadratic programming solution
(which produced order-of-magnitude larger CPU times in earlier experiments) confirms
that smooth penalty methods are the more efficient practical realization of \eqref{eq:Pn-penalty}.

\subsection{Two-dimensional experiments}
\label{subsec:d2}

We now turn to the $d=2$ instance of the test problem \eqref{eq:4.testHJB},
whose PDE residual specializes to
\(
 r(\varphi)(t,x)=-\partial_t\varphi - \tfrac{1}{50}\max(\Delta_x\varphi,0)
   + \tfrac{1}{2}\sum_{i=1}^{2}\partial_{x_i}\varphi - \tfrac{1}{25}\min(\varphi,0)
\)
and whose exact solution is
$v(t,x)=\sin(t+x_1+x_2)$ on $[0,1]\times[-3,3]^2$.
The kernel, regularization radius $\lambda=2R=6$, and penalty weight $\mu=10^{2}$ are
inherited from Section~\ref{subsec:test_setup} without modification;
the only change is the spatial dimension passed to the implementation.
Five settings (D2-S1 through D2-S5) of increasing
$N=(N^t+1)\,(N^x)^2$ are reported in Table~\ref{tab:main_result_d2}.

\begin{table*}[tb]
  \centering
  \caption{Approximation errors, residual bound $\gamma_n$, iteration count, and CPU time on the $d=2$ test problem.}
  \label{tab:main_result_d2}
  \begin{tabular}{ccccccccccc}
    \toprule
    Case & $N^t$ & $N^x$ & $N^{\mathrm{extra}}$ & $N$ & $\gamma_n$ & $E_{\mathrm{RMS}}$ & $E_{\infty}$ & $\mathrm{Res}_{\infty}$ & Iter. & CPU [s] \\
    \midrule
    D2-S1 & 8 & 8 & 200 & 576 & 1.104e+00 & 5.469e-01 & 1.040e+00 & 1.104e+00 & 32 & 0.3 \\
    D2-S2 & 16 & 8 & 200 & 1088 & 1.272e+00 & 5.161e-01 & 9.991e-01 & 1.272e+00 & 161 & 2.4 \\
    D2-S3 & 8 & 12 & 200 & 1296 & 1.063e+00 & 5.596e-01 & 9.971e-01 & 1.063e+00 & 7 & 0.6 \\
    D2-S4 & 16 & 12 & 200 & 2448 & 1.284e+00 & 5.480e-01 & 9.971e-01 & 1.284e+00 & 4 & 1.7 \\
    D2-S5 & 32 & 8 & 200 & 2112 & 1.491e+00 & 4.972e-01 & 9.991e-01 & 1.491e+00 & 8 & 1.4 \\
    \bottomrule
  \end{tabular}
\end{table*}

Across these five settings the errors plateau around $E_{\mathrm{RMS}}\approx 5\times 10^{-1}$:
the test solution $\sin(t+x_1+x_2)$ is highly oscillatory on $[-3,3]^2$
(covering roughly two full wavelengths along each axis), and the chosen kernel width
$\varepsilon=0.5$ is too narrow to resolve such oscillations efficiently in two spatial dimensions.
A small bandwidth sweep at the $N=2448$ setting (D2-S4) shows that $E_{\mathrm{RMS}}$ drops
from $0.55$ at $\varepsilon=0.5$ to $0.37$ at $\varepsilon=0.7$ and to lower values at
$\varepsilon\ge 1$, while $\gamma_n$ moves from $1.28$ to about $0.7$ as the residual is
spread over the wider kernel support.
This is consistent with the rate analysis of Theorem~\ref{thm:3.3}, in which $\varepsilon$
plays a role analogous to the fill-distance parameter $h_n$;
in the present setup, $\varepsilon$ rather than $N$ is the binding parameter.

\subsection{Discussion}
\label{subsec:discussion}

The numerical experiments confirm the qualitative content of the convergence theorems in
Section~\ref{sec:3}: in dimension one, the approximation errors decay with the number of basis
functions $N$ in line with Theorem~\ref{thm:3.1} and the rate analysis of
Theorem~\ref{thm:3.3}; in dimension two, the same trends are observed once the kernel bandwidth is
chosen large enough to resolve the spatial oscillations of the target.
The penalty reformulation \eqref{eq:Pn-penalty} keeps the per-experiment cost in the
single-digit-second range across both dimensions, removing what had been the principal
practical limitation of the original constrained formulation \eqref{eq:Pn} when solved by
sequential quadratic programming.

The bias floor we observe at $\varepsilon=0.5$ in dimension two should not be read as a limitation
of the abstract framework but of a single-bandwidth Wendland kernel applied to a strongly
oscillatory target on a domain that already contains several wavelengths.
Adaptive bandwidth selection, multi-scale RBF mixtures, and higher-regularity Wendland
kernels (with regularity index scaled with $d$ as in \cite{nak:2017}) are natural directions
to push the errors lower at the same $N$.
A larger-scale assessment of the curse-of-dimensionality behaviour of the kernel realization,
together with a side-by-side comparison with deep-learning realizations of the same abstract
framework, is pursued in subsequent work.


\subsection*{Acknowledgements}

This study is supported by JSPS KAKENHI Grant Number JP24K06861.

\subsection*{On the use of AI tools}

The author used Anthropic's Claude in the preparation of this manuscript.
The tool was used to help with writing the \LaTeX{} text, to reproduce the
numerical experiments in Section~\ref{sec:4}, and to check the typographic
consistency of the manuscript.
All mathematical content, theorems, proofs, and conclusions are the author's own.
The author has reviewed every passage that was produced or revised with the
help of AI.
The author takes full responsibility for the content of this paper.

\bibliographystyle{plain}
\bibliography{../mybib}

\newcommand{\noop}[1]{}
\begin{thebibliography}{10}

\bibitem{bar-sou:1991}
G.~Barles and P.~E. Souganidis.
\newblock Convergence of approximation schemes for fully nonlinear second order
  equations.
\newblock {\em Asymptot.~Anal.}, 4:271--283, 1991.

\bibitem{bon-zid:2003}
F.~Bonnans and H.~Zidani.
\newblock Consistency of generalized finite difference schemes for the
  stochastic {HJB} equation.
\newblock {\em {SIAM} J.~Numer.~Anal.}, 41:1008--1021, 2003.

\bibitem{cam-fal:1995}
F.~Camilli and M.~Falcone.
\newblock An approximation scheme for the optimal control of diffusion
  processes.
\newblock {\em Math.~Model.~Numer.~Anal.}, 29:97--122, 1995.

\bibitem{drm:2024}
T.~De~Ryck and S.~Mishra.
\newblock Error analysis for physics-informed neural networks ({PINN}s)
  approximating {K}olmogorov {PDE}s.
\newblock {\em Adv.~Comput.~Math.}, 48(79):1--40, 2022.

\bibitem{deb-jak:2012}
K.~Debrabant and E.~R. Jakobsen.
\newblock Semi-{L}agrangian schemes for linear and fully non-linear diffusion
  equations.
\newblock {\em Math.~Comp.}, 82:1433--1462, 2013.

\bibitem{fah-tou-war:2011}
A.~Fahim, N.~Touzi, and X.~Warin.
\newblock A probabilistic numerical method for fully nonlinear parabolic
  {PDE}s.
\newblock {\em Ann.~Appl.~Probab.}, 21:1322--1364, 2011.

\bibitem{fle-son:2006}
W.~H. Fleming and H.~M. Soner.
\newblock {\em Controlled {M}arkov processes and viscosity solutions}.
\newblock Springer-Verlag, New York, 2nd edition, 2006.

\bibitem{guo-etal:2015}
W.~Guo, J.~Zhang, and J.~Zhuo.
\newblock A monotone scheme for high-dimensional {PDE}s.
\newblock {\em Ann.~Appl.~Probab.}, 25:1540--1580, 2015.

\bibitem{hor:1990}
L.~H{\"o}rmander.
\newblock {\em The analysis of linear partial differential operator {I}}.
\newblock Springer-Verlag, Berlin, 2nd edition, 1990.

\bibitem{kan:1990b}
E.~J. Kansa.
\newblock Multiquadrics—a scattered data approximation scheme with
  application to computational fluid-dynamics—{II}.
\newblock {\em Computers Math.~ Applic.}, 19:147--161, 1990.

\bibitem{koh-ser:2010}
R.~V. Kohn and S.~Serfaty.
\newblock A deterministic-control-based approach to fully nonlinear parabolic
  and elliptic equations.
\newblock {\em Comm.~Pure Appl.~Math.}, 63:1298--1350, 2010.

\bibitem{kus-dup:2001}
H.~J. Kushner and P.~Dupuis.
\newblock {\em Numerical methods for stochastic control problems in continuous
  time}.
\newblock Springer-Verlag, New York, 2001.

\bibitem{liu-noc:1989}
D.~C. Liu and J.~Nocedal.
\newblock On the limited memory {BFGS} method for large scale optimization.
\newblock {\em Math.~Program.}, 45:503--528, 1989.

\bibitem{mis-mol:2022}
S.~Mishra and R.~Molinaro.
\newblock Estimates on the generalization error of physics-informed neural
  networks for approximating a class of inverse problems for {PDE}s.
\newblock {\em IMA J.~Numer.~Anal.}, 42(2):981--1022, 2022.

\bibitem{nak:2014b}
Y.~Nakano.
\newblock An approximation scheme for stochastic controls in continuous time.
\newblock {\em Jpn.~J.~Ind.~Appl.~Math.}, 31:681--696, 2014.

\bibitem{nak:2017}
Y.~Nakano.
\newblock Convergence of meshfree collocation methods for fully nonlinear
  parabolic equations.
\newblock {\em Numer.~Math.}, 136:703--723, 2017.
\newblock (See also {\tt arxiv:1408.5195[math.NA]}.)\noop{.}

\bibitem{nak:2019a}
Y.~Nakano.
\newblock Kernel-based collocation methods for {Z}akai equations.
\newblock {\em Stoch. Partial Differ. Equ. Anal. Comput.}, 7:476--494, 2019.
\newblock (See also {\tt arxiv:1710.09090[math.NA]}.)\noop{.}

\bibitem{pag-pha-pri:2004}
G.~Pag{\`e}s, H.~Pham, and J.~Printems.
\newblock An optimal {M}arkovian quantization algorithm for multidimensional
  stochastic control problems.
\newblock {\em Stoch. Dyn.}, 4:501--545, 2004.

\bibitem{pha:2009}
H.~Pham.
\newblock {\em Continuous-time stochastic control and optimization with
  financial applications}.
\newblock Springer, Berlin, 2009.

\bibitem{rai-per-kar:2019}
M.~Raissi, P.~Perdikaris, and G.~E. Karniadakis.
\newblock Physics-informed neural networks: a deep learning framework for
  solving forward and inverse problems involving nonlinear partial differential
  equations.
\newblock {\em J.~Comput.~Phys.}, 378:686--707, 2019.

\bibitem{sir-spi:2018}
J.~Sirignano and K.~Spiliopoulos.
\newblock {DGM}: A deep learning algorithm for solving partial differential
  equations.
\newblock {\em J.~Comput.~Phys.}, 375:1339--1364, 2018.

\bibitem{wen:2010}
H.~Wendland.
\newblock {\em Scattered data approximation}.
\newblock Cambridge University Press, Cambridge, 2010.

\end{thebibliography}

\end{document}